\documentclass[11pt,psamsfonts, two-sided]{amsart}
 
\usepackage[headinclude]{typearea}
\usepackage{amsfonts, amssymb, enumerate, mathpazo} 
\usepackage[matrix,curve]{xypic}
\usepackage[english]{babel}
\usepackage{hyperref}

\numberwithin{equation}{section}  
\areaset[1cm]{5.8in}{8.5in}

\newcommand{\ZZ}{\mathbb{Z}}
\newcommand{\CC}{\mathbb{C}}
\newcommand{\RR}{\mathbb{R}}
\newcommand{\QQ}{\mathbb{Q}}
\newcommand{\NN}{\mathbb{N}}

\newcommand{\lieg}{\mathfrak{g}}
\newcommand{\lieh}{\mathfrak{h}}

\newcommand{\FC}{\mathcal{F}}
\newcommand{\UC}{\mathcal{U}}

\newcommand{\NC}{\mathcal{N}}
\newcommand{\AC}{\mathcal{A}}
\newcommand{\BC}{\mathcal{B}}

\newcommand{\CM}{\mathcal{CM}}
\newcommand{\Ccal}{\mathcal{C}}
\newcommand{\myrtimes}{\!\rtimes\!}

\DeclareMathOperator{\tor}{Tor}
\DeclareMathOperator{\ext}{Ext}
\DeclareMathOperator{\id}{id}

\DeclareMathOperator{\im}{im}

\DeclareMathOperator{\closure}{clos}

\DeclareMathOperator{\res}{res}

\DeclareMathOperator{\boldP}{\mathbf{P}}
\DeclareMathOperator{\boldT}{\mathbf{T}}
\DeclareMathOperator{\trace}{tr}

\DeclareMathOperator{\iu}{\underline{I}}
\DeclareMathOperator{\io}{\overline{I}}
\DeclareMathOperator{\cd}{cd}
\DeclareMathOperator{\hd}{hd}
\DeclareMathOperator{\neurtimes}{\bar{\rtimes}}
\DeclareMathOperator{\dom}{dom}
\DeclareMathOperator{\ran}{ran}
\DeclareMathOperator{\redhomologie}{\overline{H}}
\DeclareMathOperator{\homologie}{H}

\newcommand{\betti}{b^{(2)}}

\newcommand{\scalar}[2]{{\langle #1,#2\rangle}}

\newcommand{\lxb}{{L^\infty (X)}}

\theoremstyle{definition}
\newtheorem{defi}{Definition}[section]

\newtheorem{rem}[defi]{Remark}

\theoremstyle{plain}
\newtheorem{lem}[defi]{Lemma}
\newtheorem{thm}[defi]{Theorem}
\newtheorem{cor}[defi]{Corollary}

\begin{document}
\title{Homological Invariants and Quasi-Isometry} 
\author{Roman Sauer}
\email{roman.sauer@uni-muenster.de}
\urladdr{www.romansauer.de}       
\address{FB Mathematik, Universit\"at M\"unster, Einsteinstr. 62,
  48149 M\"unster\\Germany} 

\begin{abstract} 
Building upon work of Y.~Shalom 
we give a homological-algebra flavored definition of an induction 
map in group homology associated to a topological coupling. 
As an application we obtain that the cohomological dimension $\cd_R$
over a commutative 
ring $R$ satisfies the inequality $\cd_R(\Lambda)\le\cd_R(\Gamma)$ 
if $\Lambda$ embeds uniformly into $\Gamma$ and
$\cd_R(\Lambda)<\infty$ holds. Another consequence of our results is 
that the Hirsch ranks of quasi-isometric solvable groups coincide. 
Further, it is shown that the real cohomology rings of quasi-isometric 
nilpotent groups are isomorphic as graded rings. On the analytic
side, we apply the induction technique to Novikov-Shubin invariants 
of amenable groups, which can be seen as homological invariants, 
and show their invariance under quasi-isometry. 
\end{abstract}

\maketitle

\section{Introduction and Statement of Results}

What is the relation between the (co)homology groups of two given 
quasi-iso\-metric groups? This question is as natural as challenging
since a quasi-isometry, being a purely geometric notion, provides 
no obvious way to produce a reasonable map in group cohomology. 
Recently, Yehuda Shalom~\cite{shalom(2003)} introduced 
a whole new circle of ideas and techniques into the study of the 
geometry of amenable groups. These techniques involve representation
theory and cohomology. The starting point 
is an induction map in group cohomology (for both ordinary and
reduced) associated to a topological 
coupling of quasi-isometric groups (defined
in~\ref{dynamic-criterion}). The existence of such a topological 
coupling is equivalent to being quasi-isometric due to 
Gromov's dynamic characterization of
quasi-isometry~\cite[$0.2.C_2'$]{gromov(1993)}. Restricting to 
ordinary cohomology, Shalom's induction map is available for the 
more general situation of uniform embeddings. 
It is defined explicitly in terms of the standard homogeneous
resolution.  \\[\smallskipamount]
We provide a different and more abstract approach to the
induction which works for ordinary cohomology \textit{and} 
homology, is natural 
in the coefficients and compatible with cup- and cap-products. 
In doing so we build heavily upon the ideas of Shalom. 
Let us discuss now the main results of the paper. 
We first recall 
the notion of uniform embedding, 
a notion encompassing 
subgroup inclusions and quasi-isometric embeddings. 

\pagebreak
\begin{defi}\label{definition-of-uniform-embedding} 
Let $\Lambda,\Gamma$ be discrete countable groups. 
\begin{enumerate}[(i)]
\item A map $\phi:\Lambda\rightarrow\Gamma$ is called a \textbf{uniform 
embedding} if for every sequence of pairs
$(\alpha_i,\beta_i)\in\Lambda\times\Lambda$ one has: 
\[\alpha_i^{-1}\beta_i\rightarrow\infty\text{ in
  $\Lambda$}\iff\phi(\alpha_i)^{-1}\phi(\beta_i)\rightarrow\infty\text{
  in $\Gamma$.}\]
Here $\rightarrow\infty$ means eventually leaving every finite
  subset. 
\item The groups $\Lambda,\Gamma$ are called \textbf{quasi-isometric} if there 
exists a uniform embedding $\phi:\Lambda\rightarrow\Gamma$ and finite 
subset $C\subset\Gamma$ such that $\phi(\Lambda)\cdot C=\Gamma$. 
\end{enumerate}
\end{defi}

That the latter definition of quasi-isometry is equivalent to 
the usual one (see definition~\ref{definition-of-quasi-isometry}) 
is part of theorem~\ref{dynamic-criterion}. Now suppose that 
$\Lambda$ uniformly embeds into $\Gamma$, and let $R$ be a
commutative ring.  
What one finally gets (theorem~\ref{theorem-induction-in-homology}),  
taking a long route via Gromov's dynamic
criterion~\ref{dynamic-criterion}, is the following: There 
is a compact topological 
space $Y$ with a continuous $\Lambda$-action, a 
functor 
$\io:\{\text{$R\Lambda$-modules}\}\rightarrow\{\text{$R\Gamma$-modules}\}$ 
and a homomorphism, called the \textbf{induction}, in cohomology 
\begin{equation}\label{eq:statement-induction}
I^n:\xymatrix@1{\homologie^n(\Lambda,M)\ar[r]&
\homologie^n(\Lambda,\FC(Y;R)\otimes_R M)\ar[r]^-{\cong} &\homologie^n(\Gamma,\io(M))}
\end{equation}
for every $R\Lambda$-module $M$. Here $\FC(Y;R)$ is the ring  
of functions $Y\rightarrow R$ 
with the property that the preimage of any $r\in R$ is 
open and closed; it carries a natural $\Lambda$-action. 
We consider $\FC(Y;R)\otimes_R M$ as an $R\Lambda$-module by the 
diagonal $\Lambda$-action. The first map
in~(\ref{eq:statement-induction}) is induced by the inclusion 
$M\hookrightarrow\FC(Y;R)\otimes_R M,~m\mapsto \id_Y\otimes m$. 
The second map in~(\ref{eq:statement-induction}) is always an
isomorphism. 
If we can prove that  
the first map and hence $I^n$ are 
injective under certain assumptions, then we get the 
estimate $\cd_R(\Lambda)\le\cd_R(\Gamma)$ for the cohomological
dimensions over $R$. So we do not need to care so much about the 
actual definitions of $\io(M)$ and the second map
in~(\ref{eq:statement-induction}). The theorem we obtain by analyzing
the first map is the 
following (shown in section~\ref{sec:quasi-isometry and homological
  dimension}). 

\begin{thm}\label{statement-theorem estimate with homological
    dimension} 
Let $R$ be a commutative ring, and suppose $\Lambda$ embeds uniformly 
into $\Gamma$ where $\Lambda$ and $\Gamma$ are discrete, countable groups. 
Then the following two statements hold. 
\begin{enumerate}[(i)]
\item If $\cd_R(\Lambda)$ is finite, then we have
  $\cd_R(\Lambda)\le\cd_R(\Gamma)$. 
\item If $\Lambda$ is amenable and $\QQ\subset R$, then we have 
  $\cd_R(\Lambda)\le\cd_R(\Gamma)$. 
\end{enumerate}
Furthermore, $(i),(ii)$ hold true if $\cd_R$ is replaced by the
homological dimension $\hd_R$. 
\end{thm}

Here statement (ii) for the cohomological dimension is already 
proved in~\cite[theorem 1.5]{shalom(2003)} and was 
conjectured for 
the homological dimension in~\cite[section 6.4]{shalom(2003)}. An important 
point is that we can deal with non-amenable groups by imposing 
a finiteness condition. The theorem above also generalizes a result of 
Gersten~\cite{gersten(1993)}.  
\smallskip\\
By a result of Stammbach~\cite{stammbach(1970)}, 
the rational homological dimension of a 
solvable group equals its Hirsch number. Hence we obtain the 
following corollary which was known before only under 
additional finiteness conditions on the groups
(see~\cite{bridson(1996)}, \cite{shalom(2003)}).  

\begin{cor}\label{statement-theorem Hirsch rank} 
Let $\Gamma$ be a solvable group, and let $\Lambda$ be a solvable 
group quasi-isometric to $\Gamma$. Then the Hirsch ranks of $\Gamma$
and $\Lambda$ coincide. 
\end{cor}

Recall that the Hirsch rank $h(\Gamma)$ 
of a solvable group $\Gamma$ is defined 
as the sum 
\[h(\Gamma)=\sum_{i\ge
  0}\dim_\QQ\bigl(\Gamma^{i}/\Gamma^{i+1}\otimes_\ZZ\QQ\bigr),\]
 where $\Gamma^{(i)}$ is the $i$-th term in the derived series of $\Gamma$. 
An interesting application of~\ref{statement-theorem estimate with homological
    dimension}~(i) was pointed out to me 
by Shalom. Let be $\Gamma,\Lambda$ two 
non-uniform arithmetic lattices in the same connected semisimple Lie group 
$G$ with finite center. 
Conjecturally~\cite[section 6.1]{shalom(2003)},  
$\Lambda$ uniformly embeds into $\Gamma$ if and only if 
$\text{$\QQ$-rank}~\Lambda\ge\text{$\QQ$-rank}~\Gamma$, with equality 
if and only if $\Lambda,\Gamma$ are commensurable. We remark that 
if $\Gamma$ is uniform, i.e.$~\text{$\QQ$-rank}~\Gamma=0$, then every 
discrete subgroup of $G$ embeds uniformly into $\Gamma$. 
Then~\ref{statement-theorem estimate with homological 
dimension}~(i) 
together with a theorem of Borel and Serre~\cite{borelserre(1973)}, 
which expresses $\cd_\QQ(\Gamma)$ of a lattice $\Gamma$ as the 
difference of the dimension of the associated symmetric space with
the $\QQ$-rank of $\Gamma$, imply the "only if"-statement:  

\begin{cor}\label{statement-cor uniform embedding and lattices} 
Let $\Gamma,\Lambda$ be arithmetic lattices in the same connected 
semisimple Lie group with finite center. 
If $\Lambda$ uniformly embeds into $\Gamma$ then 
$\text{$\QQ$-rank}~\Lambda\ge\text{$\QQ$-rank}~\Gamma$ holds. 
\end{cor}

We remark that irreducible lattices in semisimple 
Lie groups of $\text{$\RR$-rank}\ge 2$ are arithmetic by Margulis'
arithmeticity theorem. \\[\smallskipamount]
An interesting generalization of Shalom's theorem 
saying that the Betti numbers of quasi-isometric 
nilpotent groups coincide~\cite[theorem 1.2]{shalom(2003)} 
is proved in 
section~\ref{subsec:proof of theorem on cohomology ring of nilpotent
  groups}:  

\begin{thm}\label{thm-statement: betti numbers of nilpotent groups} 
If $\Gamma$ and $\Lambda$ are quasi-isometric nilpotent groups, 
then the real cohomology rings $\homologie^\ast(\Gamma,\RR)$ and 
$\homologie^\ast(\Lambda,\RR)$ are isomorphic as graded rings. 
\end{thm}

By a theorem of Malcev any finitely generated torsion-free nilpotent 
group $\Gamma$ is discretely and cocompactly embedded in a unique 
simply connected nilpotent Lie group $G$, the so-called 
\textit{(real) Malcev completion} of $\Gamma$. Thus $\Gamma$ has an 
associated real Lie algebra $\mathfrak{g}$. By~\cite[theorem
1]{nomizu(1954)} the cohomology algebras of $\mathfrak{g}$ and 
$\Gamma$ are isomorphic: 
\[\homologie^\ast(\Gamma,\RR)\cong
\homologie^\ast(\mathfrak{g},\RR).\]
It is a long standing question whether the real Malcev completions 
of quasi-isometric nilpotent groups are isomorphic. Note that 
a positive answer would imply the preceding theorem. We remark 
that the \textit{graded} Lie algebra associated to $\mathfrak{g}$ is 
a quasi-isometry invariant of $\Gamma$ 
by a deep theorem of Pierre Pansu~\cite{pansu(1989)}.\\[\smallskipamount]

On the analytic side, we apply our methods to \textbf{Novikov-Shubin 
invariants} of ame\-nable groups, i.e.~of the classifying spaces 
of these groups. The $i$-th Novikov-Shubin 
invariant, defined by Novikov and Shubin in the 
80's~\cite{Novikov-Shubin(1986a)},~\cite{Novikov-Shubin(1986b)}, 
can be seen as a kind of secondary information 
associated to the $i$-th $L^2$-Betti number. 
In the Riemannian setting, 
the $i$-th Novikov-Shubin invariant $\alpha_i(\tilde{M})$ 
of the universal covering $\tilde{M}$ 
of a compact Riemannian manifold $M$ quantifies the speed of the convergence 
of the limit 
\[\betti_i(\tilde{M})=\lim_{t\rightarrow\infty}\underbrace{\int_\mathcal{F}
  \trace_\RR\left ( e^{-t\Delta_i}(x,x)\right)dvol_x}_{=\theta_i(t)}.\] 
Here $\FC$ is a fundamental domain for the $\pi_1(M)$-action on 
$\tilde{M}$, $\Delta_i$ is the Laplacian on the $i$-forms on $\tilde{M}$, 
and $\betti_i(\tilde{M})$ denotes the $i$-th $L^2$-Betti number of
  $\tilde{M}$. For instance, if the integral $\theta_i(t)$ has an
  asymptotic 
  behavior like $t^{-p}+\betti_i(\tilde{M})$ for
  $t\rightarrow\infty$, then $\alpha_i(\tilde{M})$ would be equal to
  $p$. 
As for $\betti_i(X)$, there is a notion of the Novikov-Shubin
invariant 
$\alpha_i(X)$ of a finite type (i.e.~finitely many $\Gamma$-cells in 
each dimension) free $\Gamma$-CW complex $X$, 
which coincides 
with the heat kernel definition in the case of universal coverings of 
compact Riemannian 
manifolds. For more information see~\cite[chapter 2]{Lueck(2002)}. 
The Novikov-Shubin invariants $\alpha_i(\Gamma),~i\ge 1$,  
of a group $\Gamma$ are 
defined as the Novikov-Shubin invariants of the classifying 
space $E\Gamma$ provided $E\Gamma$ admits a model of finite type. 
Their relation to the geometry of groups is already indicated 
by the value of $\alpha_1(\Gamma)$ 
for finitely generated $\Gamma$ 
(see the computation 
in~\cite[proposition 3.2]{Lueck-Reich-Schick(1999)}, 
based on results of Varopoulos). 
\[\alpha_1(\Gamma)=\begin{cases} \infty^+ &\text{ if $\Gamma$ is finite or non-amenable,}\\
                       n &\text{ if $\Gamma$ has polynomial growth of
                       degree $n$,}\\
                       \infty   &\text{ otherwise.}
                     \end{cases}\]
We briefly dwell on the definition of $\alpha_i(\Gamma)$ we 
will actually work with. This definition of $\alpha_i(\Gamma)$, 
developed in~\cite{Lueck-Reich-Schick(1999)}, interprets 
$\alpha_i(\Gamma)$ as an \textit{invariant of the group homology} of 
$\Gamma$ with coefficients in the group von Neumann algebra 
$\NC(\Gamma)$, and is available for any group, 
not only for those with a finite type classifying space. 
The point of view is similar as in the algebraic 
definition of $L^2$-Betti numbers by Wolfgang 
L\"uck~\cite{Lueck(1998a)},~\cite{Lueck(1998b)}, and so is 
the motivation. Here is a typical situation: 
A group $\Gamma$ to which 
the original definition of $\alpha_i(\Gamma)$ applies 
could have a normal subgroup $\Lambda$ to which it does not apply, 
but for which we know the value or an estimate of 
$\alpha_i(\Lambda)$ (extended definition). Then this information tells 
us something about $\alpha_i(\Gamma)$ by using the Hochschild-Serre 
spectral sequence (cf.~\cite{Lueck-Reich-Schick(1999)}). 

However, as 
opposed to $L^2$-Betti numbers, this extension to all groups 
does not preserve all the properties of $\alpha_i$ one wants to
have; the maximal subclass $\CM$ of amenable groups for which we get 
a reasonable notion of $\alpha_i$ is described in
definition~\ref{definition of the class CM}. The class $\CM$ contains
inter alia all amenable groups of type $FP_\infty$ over $\CC$, hence 
including the class of amenable groups for which $\alpha_i$ was  
originally defined, and is closed under quasi-isometry. 

\begin{thm}\label{statement-theorem qi invariance of novikov-shubin
    invariants} 
Let be $\Gamma\in\CM$. If $\Lambda$ is quasi-isometric to $\Gamma$, then 
$\alpha_i(\Gamma)=\alpha_i(\Lambda)$ for $i\ge 1$. 
\end{thm}

According to~\cite[8.A6]{gromov(1993)} it is not 
unreasonable to expect that the quasi-isometry 
invariance of $\alpha_i$ holds true for, at least,  
all the groups to which the classical definition of 
Novikov-Shubin invariants applies. \\[\medskipamount]

\noindent\textbf{Acknowledgments.} My gratitude goes to 
Yehuda Shalom for encouragement and his 
interest in this work. Further, I thank him for pointing out 
some inaccuracies and giving hints for improvement. 

\pagebreak
\tableofcontents

\section{Quasi-Isometry and Topological
  Couplings}\label{sec:quasi-isometry-and-topological-couplings}      

This technical section lays the basis for the whole paper. 
As in Shalom's work on the geometry of 
amenable groups our starting point is 
the dynamic viewpoint on 
quasi-isometry by Gromov~\cite[$0.2.C_2'$]{gromov(1993)}. Accordingly  
the existence of a topological coupling of groups $\Lambda,\Gamma$ (see
theorem~\ref{dynamic-criterion}) is a characterizing property of 
$\Lambda,\Gamma$ being quasi-isometric. 
In the measurable setting, it is well known that a measure equivalence 
of groups gives rise to a weak orbit equivalence. In our situation we obtain 
a sort of topological version of weak orbit equivalence (see
subsection~\ref{subsec:uniform orbit equivalence}) which induces 
an isomorphism between certain transformation groupoids of $\Lambda$ 
and $\Gamma$ (lemma~\ref{isomorphism-of-transformation-groupoids}). 
To such a transformation groupoid of a group we associate 
a ring (see subsection~\ref{subsec:algebraic objects associated to group
  actions}) containing the group ring up to restriction to an
idempotent. It is crucial that we have a good control over 
the passage from the 
group ring of $\Lambda$ resp.~$\Gamma$ to its groupoid ring
(see e.g.~lemma~\ref{flatness-of-crossed-product-ring}), and that the 
groupoid rings of $\Lambda,\Gamma$ are isomorphic. Hence we 
can compare $\Lambda$ and $\Gamma$ algebraically.

\subsection{Topological and
  couplings}\label{subsec:topological-and-uniform-measure-couplings}  
Let us recall the
standard definition of a quasi-isometry between finitely 
generated groups. 

\begin{defi}\label{definition-of-quasi-isometry} 
Let $\Gamma$, $\Lambda$ be groups generated by the finite symmetric 
sets $S_\Gamma$, $S_\Lambda$ and equipped with the corresponding 
word metrics $d_\Gamma$, $d_\Lambda$ on $\Gamma,\Lambda$. A map 
$\phi:\Lambda\rightarrow\Gamma$ is called a \textbf{quasi-isometric embedding}
if there are constants $\alpha\ge 1$, 
$C\ge 0$, such that for all $\lambda_1, \lambda_2\in\Lambda$ one has 
\[\alpha^{-1}d_\Lambda(\lambda_1,\lambda_2)-C\le
d_\Gamma(\phi(\lambda_1),\phi(\lambda_2))\le\alpha
d_\Lambda(\lambda_1,\lambda_2)+C.\]
If, in addition, any $\gamma\in\Gamma$ lies within distance $\le D$, for
a constant $D\ge 0$, from the image $\phi(\Lambda)$, then $\phi$ is 
called a \textbf{quasi-isometry}. The groups $\Gamma,\Lambda$ are called 
quasi-isometric if there exists a quasi-isometry
$\phi:\Lambda\rightarrow\Gamma$. 
\end{defi}
The following theorem is essentially Gromov's dynamic
criterion. See~\cite[$0.2.C_2'$]{gromov(1993)} and especially~\cite[theorem
2.1.2]{shalom(2003)}. 

\begin{thm}\label{dynamic-criterion}
For countable groups $\Lambda,\Gamma$ consider the following 
statements. 
\begin{enumerate}[(i)]
\item There exists a uniform embedding
  $\phi:\Lambda\rightarrow\Gamma$. 
\item There exists a locally compact space $X$ on which both $\Lambda$
  and $\Gamma$ act continuously, freely and properly such that the two actions 
commute and the $\Gamma$-action is cocompact. Further, there 
exists fundamental domains $X_\Lambda$ and $X_\Gamma$ such that 
$X_\Gamma$ is compact-open and $X_\Lambda$ is closed-open. 
The space $X$ is called a \textbf{topological coupling} of $\Lambda$ 
and $\Gamma$. 
\item There exists $\phi$ as in (i) and a finite subset
  $C\subset\Gamma$ such that $\phi(\Lambda)\cdot C=\Gamma$. 
\item There exists $X$ as in (ii) but with $X_\Lambda$ being 
compact, i.e.~the $\Lambda$-action is also cocompact. 
\end{enumerate}
Then (i) is equivalent to (ii) and (iii) is equivalent to (iv). 
Furthermore, if $\Lambda,\Gamma$ are finitely generated, 
any $\phi$ as in (iii) is a 
quasi-isometry of $\Lambda$ with $\Gamma$ and (iii) or (iv) are 
equivalent to $\Lambda$ and $\Gamma$ being quasi-isometric in 
the sense of~\ref{definition-of-quasi-isometry}. 
\end{thm}

\begin{rem}\label{shaloms-topological-coupling-version-with-inclusion}
This is the formulation in~\cite[theorem 2.1.2]{shalom(2003)} except 
that there the freeness of the actions in (ii) is not
demanded. Instead, it is shown that if one replaces $\Gamma$ with a
direct product $\Gamma\times F$ for some finite group $F$, then a space 
$X$ can be found with the same properties as (ii) and the additional 
property that $X_\Gamma\subset X_\Lambda$. For us it will be important 
to have the version above since we want to prove 
theorems about the cohomological dimension $\cd_R$ over an arbitrary ring $R$, 
so we are not allowed to replace $\Gamma$ by the product with a 
finite group without changing $\cd_R$. 
\end{rem}

\begin{proof} 
We only indicate how 
to modify the construction in~\cite[proof of 2.1.2]{shalom(2003)} so that one 
gets a topological coupling as above. 
Let $\phi:\Lambda\rightarrow\Gamma$ be a uniform embedding. Then there 
is a finite subset $Q\subset\Lambda$ such that if
$\phi(\lambda_1)=\phi(\lambda_2)$ then $\lambda_2^{-1}\lambda_1\in
Q$. If $F$ is a finite group having more elements than $Q$, then it is 
easy to see that there is an injective uniform embedding
$\phi':\Lambda\rightarrow\Gamma\times F$ such that $p_\Gamma\circ\phi'=\phi$ 
where $p_\Gamma:\Gamma\times F\rightarrow\Gamma$ is the 
projection onto $\Gamma$. Choose left invariant proper (i.e.~balls are
finite) metrics 
$d_\Lambda$, $d_{\Gamma\times F}$ on $\Lambda,\Gamma\times F$. Then define 
\begin{align*} 
F_1(t)&=\inf\{d_{\Gamma\times
  F}(\phi'(\lambda_1),\phi'(\lambda_2));~d_\Lambda(\lambda_1,\lambda_2)\ge
t\}\\ 
F_2(t)&=\sup\{d_{\Gamma\times
  F}(\phi'(\lambda_1),\phi'(\lambda_2));~d_\Lambda(\lambda_1,\lambda_2)\le t\}
\end{align*}
Now consider the space 
$X$ of all \textit{injective} maps $\psi:\Lambda\rightarrow\Gamma\times F$ 
satisfying the same uniform estimate as $\phi'$ does: 
\begin{equation}\label{eq:uniform-estimate}
F_1(d_\Lambda(\lambda_1,\lambda_2))\le
d_{\Gamma\times F}(\psi(\lambda_1),\psi(\lambda_2))\le
F_2(d_\Lambda(\lambda_1,\lambda_2))
\end{equation}
We equip $X$ with the pointwise-convergence topology. We let $\Gamma$ act 
on $\Gamma\times F$ by $\gamma(\gamma',m)=(\gamma\gamma',m)$. 
Thereby we obtain a right $\Lambda$- and left $\Gamma$-action on 
$X$ by $(\lambda\psi)(x)=\psi(\lambda x)$ and 
$(\gamma\psi)(x)=\gamma\psi(x)$ for $\psi\in X$. 
The $\Gamma$-action is free and 
a compact-open $\Gamma$-fundamental domain is given by 
$X_\Gamma=\{\psi\in X;~\psi(e)\in\{e\}\times F\}$. The 
$\Gamma$- and $\Lambda$-action are free and proper but the latter 
is not cocompact in general. However, if there is a finite
$C\subset\Gamma$ with $\phi(\Lambda)\cdot C=\Gamma$, then 
we add the condition that all maps $\psi$ in $X$ satisfy 
$(p_\Gamma\circ\psi)(\Lambda)\cdot C=\Gamma$. Then the compact 
subset 
\[K=\{\psi\in X;~\psi(e)\in C^{-1}\times F\}\] 
satisfies 
$K\Lambda=X$. Hence $\Lambda$ acts cocompactly in this case. 
As for the existence of a closed-open $\Lambda$-fundamental 
domain, choose an enumeration $\alpha_0=e, \alpha_1,\alpha_2,\ldots$ 
of $\Gamma\times F$. Then define 
\[E_i=\{\psi\in X;~\psi(e)=\alpha_i\}\]
and $K_i=\Lambda\cdot E_i$. Then a $\Lambda$-fundamental domain 
is given by 
\[X_\Lambda=E_0\cup\bigcup^\infty_{i=1}E_i\cap K_{i-1}^C\cap\cdots\cap
K_0^C,\] 
where $K_i^C=X-K_i$. So, if $\psi\in X_\Lambda$ and $n$ is the minimal integer 
such that $\psi$ takes the value $\alpha_n$ then $\psi(e)=\alpha_n$. 
Note that we do not necessarily have $X_\Gamma\subset X_\Lambda$. 
The proof that $X_\Lambda$ is closed-open and the 
proof of all the other topological properties of $X,X_\Lambda,X_\Gamma$ 
are exactly the same as in~\cite[proof of 2.1.2]{shalom(2003)}. 
\end{proof}
Crucial for applications to amenable group is the 
following fact (see~\cite[theorem 2.1.7]{shalom(2003)}). 

\begin{thm}\label{ergodic ume coupling out of a qi} 
Let $X$ be a topological coupling of quasi-isometric groups 
$\Gamma$, $\Lambda$ with fundamental domains $X_\Gamma,X_\Lambda$ 
as in~\ref{dynamic-criterion}. 
If $\Gamma,\Lambda$ are amenable, we can equip $X$ with a 
non-trivial, ergodic $\Gamma\times\Lambda$-invariant Borel measure
such that $X_\Gamma$ and $X_\Lambda$ have finite measure. 
\end{thm}

\subsection{A kind of topological orbit equivalence}\label{subsec:uniform orbit
  equivalence}

Let $X$ be a topological 
coupling of $\Gamma,~\Lambda$ as 
in~\ref{dynamic-criterion}. Since the actions on 
$X$ commute we obtain a \textbf{right $\Lambda$-action} on 
$\Gamma\backslash X\cong X_\Gamma$
and a \textbf{left $\Gamma$-action} on $X/\Lambda\cong X_\Lambda$. 
Note that we have natural homeomorphisms $X_\Lambda\times\Lambda\cong
X$ and $X_\Gamma\times\Gamma\cong X$ since $X_\Lambda$ and $X_\Gamma$ 
are closed and open. 
To avoid confusion, we use the dot-notation "$\gamma\cdot x$" only 
for the actions on the fundamental domains. 
We get a left $\Lambda$-action on $X_\Gamma$ simply by 
$\lambda\cdot x=x\cdot\lambda^{-1}$, which we often use 
instead of the right action for symmetry reasons. \\[\smallskipamount]
We adopt and recall the cocycle notation from~\cite[2.2]{shalom(2003)}. 
We define maps $\alpha:\Gamma\times X\rightarrow\Lambda$, 
$\beta:X\times\Lambda\rightarrow\Gamma$ by 
\begin{align} 
\alpha(\gamma,x)=\lambda &\iff (\gamma^{-1}x)\lambda\in
X_\Lambda\label{eq:cocycle1}\\
\beta(x,\lambda)=\gamma &\iff \gamma^{-1}(x\lambda)\in
X_\Gamma.\label{eq:cocycle2} 
\end{align}
Note that $\alpha,\beta$ are well defined because $X_\Lambda,X_\Gamma$
are fundamental domains. The natural $\Gamma$- and $\Lambda$-actions 
on the fundamental domains obtained from the identification as 
quotients of $X$ take the following forms. 
\begin{align} 
\gamma\cdot x &=\gamma x\alpha(\gamma^{-1},x),~~x\in
X_\Lambda\label{eq:action1} \\
x\cdot\lambda &=\beta(x,\lambda)^{-1}x\lambda,~~x\in X_\Gamma\label{eq:action2}
\end{align}
The maps $\alpha, \beta$ satisfy the following cocycle identities on
the fundamental domains. 
\begin{align}
\alpha(\gamma_1\gamma_2,x)&=\alpha(\gamma_1,x)\alpha(\gamma_2,\gamma_1^{-1}
\cdot x)&\forall\gamma_1,\gamma_2\in\Gamma,~x\in
X_\Lambda\label{eq:cocycle-identity-alpha}\\
\beta(x,\lambda_1\lambda_2)&=\beta(x,\lambda_1)\beta(x\cdot\lambda_1,\lambda_2)&\forall\lambda_1,\lambda_2\in\Lambda,~x\in   
X_\Gamma\label{eq:cocycle-identity-beta} 
\end{align}

\begin{defi}\label{definition closed-open topology} 
Note that the collection of subsets of a 
topological space, which are open and closed, forms a set algebra. We call 
a map $f:X\rightarrow Y$ between topological spaces 
\textbf{cut-and-paste continuous} if $f$ is measurable with respect 
to the set algebras of closed-open subsets. 
\end{defi}

The following lemma is a topological version of the fact from ergodic 
theory that a measure 
equivalence gives rise to a weak orbit equivalence. Unfortunately, 
we need the explicit formulation below in
lemma~\ref{isomorphism-between-groupoid-rings}, which accounts 
for some cumbersome cocycle calculations.

\begin{lem}\label{lemma-uniform-orbit-equivalence} 
Let $X$ be a topological coupling of $\Lambda, \Gamma$ 
as in theorem~\ref{dynamic-criterion} (ii) with closed-open 
fundamental domains $X_\Lambda, X_\Gamma$ where $X_\Gamma$ is 
compact. 
Let $p:X\rightarrow X_\Gamma=\Gamma\backslash X$ and $q:X\rightarrow
X_\Lambda=X/\Lambda$ be the canonical projections. 
Then there are compact-open 
subsets $A\subset X_\Gamma, B\subset X_\Lambda$ that meet every 
orbit of the action of $\Lambda$ on $X_\Gamma$ resp.~$\Gamma$ on 
$X_\Lambda$ such that the restriction $ f$ of $q$ to $A\subset
X_\Gamma$ is a cut-and-paste continuous bijection 
$ f:A\rightarrow B$ satisfying the identities 
\begin{align} 
 f(\lambda\cdot
a)&=\beta(a,\lambda^{-1})^{-1}\cdot f(a)\label{eq:equivariance-property-phi}
\\ 
 f^{-1}(\gamma\cdot
b)&=\Bigl(\alpha\bigl(1, f^{-1}(b)\bigr)\alpha\bigl(\gamma^{-1},b\bigr)\alpha
\bigl(1, f^{-1}(\gamma\cdot b)\bigr)^{-1}\Bigr)^{-1}
\cdot f^{-1}(b)
\label{eq:equivariance-property-inverse-phi}
\end{align}
for $a\in A\cap\lambda^{-1}\cdot A$, $b\in B\cap\gamma^{-1}\cdot B$. 
Furthermore, if $X$ carries a $\Gamma\times\Lambda$-invariant 
Borel measure $\mu$, then $f$ is measure-preserving with respect to 
the restrictions $\mu_{\vert A}$, $\mu_{\vert B}$. 
\end{lem}

\begin{proof} 
The image $B=q(X_\Gamma)$ is clearly compact and open. Note that under the
identification $X=X_\Lambda\times\Lambda$ there is a finite subset 
$F\subset\Lambda$ such that $X_\Gamma\subset X_\Lambda\times F$, 
in other words, $X_\Gamma$ is covered by finitely many
$\Lambda$-translates of $X_\Lambda$ since $X_\Gamma$ is compact and 
$X_\Lambda$ is open. Now it is easy to see that there is 
a closed-open subset $A\subset X_\Gamma$ such that the restriction 
of the projection 
$f=q_{\vert A}:A\rightarrow B$ is a bijection and $q(A)=q(X_\Gamma)$ holds. 
For an  
invariant measure $\mu$ on $X$ the map 
$f$ is measure-preserving since it is given by "cutting" $A$ into 
finitely many pieces and translating them by elements of 
$F\subset\Lambda$. 
From $X_\Lambda\subset\Gamma X_\Gamma$ and the $\Gamma$-equivariance of $q$ 
follow that 
$\Gamma\cdot B=X_\Lambda$, so $B$ meets every $\Gamma$-orbit. 
Further, $\Lambda\cdot A=X_\Gamma$ is obtained from $q(A)=q(X_\Gamma)$. 
\smallskip\\
Concerning the
properties~(\ref{eq:equivariance-property-phi}),~(\ref{eq:equivariance-property-inverse-phi}), we can conclude from $q$ being 
$\Gamma$-equivariant that 
\[ f(\lambda\cdot
a)= f(a\cdot\lambda^{-1})=q(\beta(a,\lambda^{-1})^{-1}a\lambda^{-1})=
\beta(a,\lambda^{-1})^{-1}\cdot
q(a\lambda^{-1})=\beta(a,\lambda^{-1})^{-1}\cdot f(a).\]
For the corresponding property of $ f^{-1}$, consider 
$b\in B$ with $\gamma\cdot b\in B$. So there are $x,w\in A$ with 
$ f(x)=b,  f(w)=\gamma\cdot f(x)=\gamma\cdot b$. 
Equation~(\ref{eq:equivariance-property-inverse-phi}) 
would follow from 
\[w=\Bigl(\alpha(1,x)\alpha(\gamma^{-1}, f(x))\alpha(1,w)^{-1}\Bigr)^{-1}x.
\] 
Note that we have 
$ f(x)=x\alpha(1,x)$ and $ f(w)=w\alpha(1,w)$ 
from which we obtain 
\begin{align*}
w= f(w)\alpha(1,w)^{-1}=\bigl(\gamma\cdot f(x)\bigr)\alpha(1,w)^{-1}&=
\gamma f(x)\alpha(\gamma^{-1}, f(x))\alpha(1,w)^{-1}\\
&=\gamma x\alpha(1,x)\alpha(\gamma^{-1}, f(x))\alpha(1,w)^{-1}\\
&=x\cdot\bigl(\alpha(1,x)\alpha(\gamma^{-1}, f(x))\alpha(1,w)^{-1}\bigr)
\end{align*} 
as desired. The last equality follows from the 
identities~(\ref{eq:cocycle2}),~(\ref{eq:action2}).  
\end{proof}

\subsection{Transformation groupoids}\label{subsec:transformation groupoids}

\begin{defi}\label{definition-action-by-pseudogroups} 
Let $Y$ be a topological space. Define $PI(Y)$ to be the set of 
partial bijections $f$ of $Y$ from a closed-open 
subset $\dom(f)$ onto 
a closed-open subset $\ran(f)$ such that $f$, $f^{-1}$ are 
cut-and-paste continuous. A 
\textbf{cut-and-paste continuous pseudo-action} of a group $\Gamma$ on 
$X$ is a map $\sigma:\Gamma\rightarrow PI(Y)$ such that 
$\sigma(1_\Gamma)=\id_Y$ and $\sigma(\gamma_1\gamma_2^{-1})=\sigma(\gamma_1)\circ\sigma(\gamma_2)^{-1}$ hold, 
where $\circ$ is the composition of partial maps. 
\end{defi}

\begin{rem}\label{typical-example-of-pseudo-action} 
The only example of a pseudo-action we will 
consider arises as the restriction of an ordinary continuous group 
action of a group $\Gamma$ on a topological 
space $Y$ to an open-closed subset $A\subset Y$. 
This pseudo-action is explicitly defined by 
$\sigma:\Gamma\rightarrow PI(A)$ sending $\gamma\in\Gamma$ 
to $\sigma(\gamma):A\cap\gamma^{-1}A\rightarrow A\cap\gamma A,~
y\mapsto\gamma y$. 
\end{rem}

\begin{defi}\label{definition-of-transformation-groupoid} 
Let $Y$ be a topological space equipped with a cut-and-paste continuous  
pseudo-action $\sigma:\Gamma\rightarrow PI(Y)$. 
Then we define the \textbf{transformation 
groupoid} $Y\rtimes\Gamma$ as 
\[Y\rtimes\Gamma=\{(y,\gamma);~y\in \dom(\sigma(\gamma)), \gamma\in\Gamma\}.\]
The unit space of this groupoid is given by $Y$, the source and 
target maps are 
$s(y,\gamma)=y$, $t(y,\gamma)=\sigma(\gamma)(y)=\gamma y$ and 
the product and the inverse are defined by 
\begin{align*}
(x,\gamma)(y,\gamma')&=(x,\gamma'\gamma)\\
(y,\gamma)^{-1}&=\bigl(\gamma y,\gamma^{-1}\bigr).
\end{align*}
\end{defi}

\begin{lem}\label{isomorphism-of-transformation-groupoids} 
We retain the notation of lemma~\ref{lemma-uniform-orbit-equivalence}. The 
map 
\[\phi:A\rtimes\Lambda\longrightarrow B\rtimes\Gamma,~(a,\lambda)\mapsto
\Bigl( f(a), \beta\bigl(a,\lambda^{-1}\bigr)^{-1}\Bigr)\]
is an isomorphism of groupoids. Its inverse is explicitly given by 
\[ \phi^{-1}(b,\gamma)=\left( f^{-1}(b),\Bigl(\alpha\bigl(1, f^{-1}(b)\bigr)
\alpha\bigl(
\gamma^{-1},b\bigr)\alpha\bigl(1, f^{-1}(\gamma\cdot b)\bigr)^{-1}\Bigr)^{-1}\right).\]
\end{lem}

\begin{proof} 
Lemma~\ref{lemma-uniform-orbit-equivalence} yields that $\phi$ is well
defined, i.e.~$ f(a)\in \dom(\sigma(\beta(a,\lambda^{-1})^{-1}))$
since $\beta(a,\lambda^{-1})^{-1}\cdot  f(a)= f(\lambda\cdot a)\in
B$. Further, $\phi$ is a groupoid morphism because of the cocycle
identity~(\ref{eq:cocycle-identity-beta}):
\begin{align*}
\phi\bigl((a,\lambda)(\lambda\cdot
a,\lambda')\bigr)&=\phi\bigl((a,\lambda'\lambda)\bigr)\\
&=\Bigl( f(a),\beta\bigl(a,\lambda^{-1}\lambda'^{-1}\bigr)^{-1}\Bigr)\\
&=\Bigl( f(a),\beta\bigl(\lambda\cdot
a,\lambda'^{-1}\bigr)^{-1}\beta\bigl(a,\lambda^{-1}\bigr)^{-1}\Bigr)\\ 
&=\phi(a,\lambda)\phi(\lambda\cdot a,\lambda')
\end{align*}
Next we only show that $\phi^{-1}\circ\phi=\id$ as 
$\phi\circ\phi^{-1}=\id$ follows along the same lines. 
For that we have to check the following cocycle identity 
\begin{equation}\label{eq:cocycle-identity-for-proving-the-inverse-phi}
\alpha(1,a)\alpha\bigl(\beta\bigl(a,\lambda^{-1}\bigr),b\bigr)\alpha(1,\lambda\cdot
a)^{-1}=\lambda^{-1}
\end{equation}
for $a\in A$ and $b= f(a)$. 
From~(\ref{eq:cocycle1}) and~(\ref{eq:action2}) we obtain 
\[\beta\bigl(a,\lambda^{-1}\bigr)^{-1}a\lambda^{-1}\alpha(1,\lambda\cdot a)\in
X_\Lambda.\] 
Note that $b= f(a)=a\alpha(1,a)\in X_\Lambda$. So we obtain
from~(\ref{eq:cocycle1}) that 
\[\beta\bigl(a,\lambda^{-1}\bigr)^{-1}a\alpha(1,a)\alpha\bigl(\beta\bigl(a,\lambda^{-1}\bigr),b\bigr)=\beta\bigl(a,\lambda^{-1}\bigr)^{-1}b\alpha\bigl(\beta
\bigl(a,\lambda^{-1}\bigr),b\bigr)\in X_\Lambda.\]
The two preceding equations and the fact that 
$X_\Lambda$ is a $\Lambda$-fundamental domain together
imply~(\ref{eq:cocycle-identity-for-proving-the-inverse-phi}). 
\end{proof}

\subsection{Some algebraic objects associated to group
  actions}\label{subsec:algebraic objects associated to group actions}

\begin{defi}\label{definition-of-groupoid-ring} 
Let $R$ be a ring, and be $Y$ a topological space equipped with a 
cut-and-paste continuous pseudo-action $\sigma$ of $\Gamma$. 
The \textbf{algebraic groupoid ring} $R(Y\rtimes\Gamma)$ 
of the transformation 
groupoid $Y\rtimes\Gamma$ is the set of all functions 
$f:Y\rtimes\Gamma\rightarrow R$ with the following properties. 
\begin{enumerate}[(i)]
\item For any $\gamma\in\Gamma$, the map 
$y\mapsto f(y,\gamma)$ from $\dom(\sigma(\gamma))$ to $R$ is
cut-and-paste continuous in the sense that the preimages of all 
$r\in R$ are open and closed. 
\item There is a finite subset $F\subset\Gamma$ such that 
$f(y,\gamma)=0$ if $\gamma\not\in F$. 
\end{enumerate}
Then $R(Y\rtimes\Gamma)$ becomes a ring by pointwise addition and 
the convolution product 
\[(fg)(x,\gamma)=\sum_{(x,\gamma')(y,\gamma'')=(x,\gamma)}
f(y,\gamma'')g(x,\gamma').\]
Note that (ii) ensures that the sum is finite. Now assume that 
$Y$ is equipped with a finite Borel measure $\mu$ which is invariant 
under the pseudo-action of $\Gamma$. Then $\CC(Y\rtimes\Gamma)$ comes 
equipped with the normalized \textbf{trace} 
\[\trace(f)=\frac{1}{\mu(Y)}\int_{Y\times\{1_\Gamma\}}fd\mu.\] 
\end{defi}

\noindent The notion \textit{trace} is justified by the following lemma. 

\begin{lem}\label{trace-property}
The functional $\trace:\CC(Y\rtimes\Gamma)\rightarrow\CC$ satisfies 
the trace property, i.e.~$\trace(fg)=\trace(gf)$ for 
$f,g\in\CC(Y\rtimes\Gamma)$. 
\end{lem}

\begin{proof}
We have 
\[\trace(fg)=\sum_{\gamma\in\Gamma}\int_{\dom(\sigma(\gamma))}f(\gamma
x,\gamma^{-1})g(x,\gamma)d\mu(x).\]
Because of the invariance of $\mu$ this equals 
\begin{align*}
\sum_{\gamma\in\Gamma}\int_{\dom(\sigma(\gamma^{-1}))}f(y,\gamma^{-1})g(\gamma^{-1}y,\gamma)d\mu(y)
&=\sum_{\gamma\in\Gamma}\int_{\dom(\sigma(\gamma))}g(\gamma
y,\gamma^{-1})f(y,\gamma)d\mu(y)\\
&=\trace(gf). 
\end{align*}
\end{proof}

\begin{defi}\label{definition-of-crossed-product-ring}
The ring of cut-and-paste continuous functions $Y\rightarrow R$ is 
denoted by $\FC(Y;R)$, which is the same as $R(Y\rtimes 1)$ in 
the notation above. 
If $\Gamma$ is acting continuously on $Y$ and hence on 
$\FC(Y;R)$ by $f^\gamma(y)=f(\gamma^{-1}y)$, one can define the 
\textbf{crossed product ring} $\FC(Y;R)\rtimes\Gamma$ as follows. 
As an abelian group $\FC(Y;R)\rtimes\Gamma$ is the free 
$\FC(Y;R)$-module $\FC(Y;R)[\Gamma]$ with basis $\Gamma$. 
Its multiplication is determined by the rule 
$f\gamma g\gamma'=fg^\gamma\gamma\gamma'$ for $f,g\in\FC(Y;R)$ and 
$\gamma,\gamma'\in\Gamma$. 
\end{defi}

Obviously, $\FC(Y;R)$ is a subring of $R(Y\rtimes\Gamma)$,
and it is clear from property (ii) in the definition of the algebraic 
groupoid ring that $R(Y\rtimes\Gamma)$ is a free $\FC(Y;R)$-module 
with basis $\Gamma$ -- like the crossed product ring. Indeed, the 
map below identifies the two rings. 

\begin{rem}\label{crossed-product-ring-equals-algebraic-groupoid-ring} 
The map defined by 
\[\FC(Y;R)\rtimes\Gamma\longrightarrow R(Y\rtimes\Gamma),~\sum_\gamma
f_\gamma\gamma\mapsto \bigl((y,\gamma)\mapsto f_\gamma(\gamma y)\bigr)\]
is a natural ring isomorphism. Moreover, if $Y$ carries a finite 
$\Gamma$-invariant Borel measure $\mu$, then 
$\trace(\sum_\gamma f_\gamma\gamma)=\mu(Y)^{-1}\int_Yf_1d\mu$ defines 
a trace on $\FC(Y;R)\rtimes\Gamma$, and the isomorphism above clearly 
preserves the traces. 
The crossed product has the advantage 
of being a universal construction, whereas
theorem~\ref{isomorphism-between-groupoid-rings} is more easily seen 
with the groupoid ring. 
In the sequel we 
\textit{identificate these two rings}. 
\end{rem}

\begin{rem}\label{module-structure-on-ring-of-functions} 
Via the natural isomorphism
$\FC(Y;R)\cong\FC(Y;R)\rtimes\Gamma\otimes_{R\Gamma}R$ the 
ring of functions $\FC(Y;R)$ becomes an
$\FC(Y;R)\rtimes\Gamma$-module. Explicitly, this module structure 
is given as follows. Let $f,g\in\FC(Y;R)$ and $\gamma\in\Gamma$. Then 
we have $(f\gamma)\cdot g= fg^\gamma$. 
\end{rem}

\begin{rem}\label{restricted-groupoid-ring}
For any ring $R$ and an idempotent $p\in R$, i.e.~$p^2=p$, the set 
$pRp$ is again a ring with unit $p$.
The characteristic function 
of a closed-open subset $A\subset Y$, usually denoted 
by $\chi_A$, is an idempotent in $\FC(Y;R)$. The action of $\Gamma$ on 
$Y$ restricts to a pseudo-action on $A$ (see 
remark~\ref{typical-example-of-pseudo-action}). Then there is an
obvious identification 
\[\chi_A\FC(Y;R)\rtimes\Gamma\chi_A=\chi_AR(Y\rtimes\Gamma)\chi_A=R(A\rtimes\Gamma
).\]  
\end{rem}

\begin{lem}\label{isomorphism-between-groupoid-rings}
We retain the notation 
of~\ref{isomorphism-of-transformation-groupoids} and
of~\ref{lemma-uniform-orbit-equivalence}. The map 
$\phi:A\rtimes\Lambda\longrightarrow B\rtimes\Gamma$ 
induces a ring isomorphism 
\[\Phi:\xymatrix@1{R(B\rtimes\Gamma)\ar[r]^\cong
  &R(A\rtimes\Lambda)},~g\mapsto g\circ\phi.\]
Now assume that the topological coupling $X$ is equipped with 
a non-trivial $\Gamma\times\Lambda$-invariant Borel measure $\mu$ 
such that $X_\Gamma,X_\Lambda$ have finite measure, and let $R=\CC$. 
Then the induced $\Gamma$- resp.~$\Lambda$-action on $X_\Lambda$
resp.~$X_\Gamma$ is $\mu$-preserving and 
$\Phi$ preserves the (normalized) traces.  
\end{lem}

\begin{proof} 
Since the fundamental domains $X_\Gamma,X_\Lambda$ are closed and open 
and the actions are continuous, the sets 
$\{x;~\alpha(\gamma,x)=\lambda\}$, $\{x;~\beta(x,\lambda)=\gamma\}$
are closed  
and open for fixed 
$\gamma\in\Gamma,\lambda\in\Lambda$. It follows from the explicit
formula in 
lemma~\ref{isomorphism-of-transformation-groupoids} 
that $g\circ\phi$ is cut-and-paste continuous for 
$g\in R(B\rtimes\Gamma)$. Similarly, $h\circ\phi^{-1}$ is
cut-and-paste continuous for $h\in R(A\rtimes\Lambda)$. 
For any compact $K\subset X$ and 
fixed $\gamma\in\Gamma$, $\lambda\in\Lambda$ the (restrictions of the) 
cocycles $\alpha:\{\gamma\}\times K\rightarrow\Lambda$, 
$\beta:K\times\{\lambda\}\rightarrow\Gamma$ have finite image since 
both fundamental domains are open. Again using the explicit formulas, 
this implies that $g\circ\phi$ and 
$h\circ\phi^{-1}$ satisfy property (ii) in 
definition~\ref{definition-of-groupoid-ring}. So $\Phi$ is 
well defined with the map induced by $\phi^{-1}$ as its inverse. 
Finally suppose $X$ comes equipped with a measure $\mu$ as described 
in the hypothesis. Since $\mu$ is $\Gamma\times\Lambda$-invariant 
it is clear from~\ref{eq:action2} 
that the induced actions on the fundamental domains 
are measure preserving. Now let $g\in
R(B\rtimes\Gamma)$. For the following computation note that 
$\beta(a,1)=1$ for $a\in A$ and that the map $f:A\rightarrow B$ 
from lemma~\ref{lemma-uniform-orbit-equivalence} 
is measure preserving, in particular $\mu(A)=\mu(B)$. Therefore we obtain 
\begin{align*}
\trace(g\circ\phi) &=\frac{1}{\mu(A)}\int_Ag(f(a),\beta(a,1)^{-1})d\mu(a)\\
&=\frac{1}{\mu(A)}\int_Ag(f(a),1)d\mu(a)\\
&=\frac{1}{\mu(B)}\int_Bg(b,1)d\mu(b)=\trace(g).
\end{align*}

\end{proof}

\begin{lem}\label{flatness-of-ring-of-functions} 
For any ring $R$ and a compact topological space $Y$ 
the ring of cut-and-paste continuous functions $\FC(Y;R)$ is a flat 
$R$-module. 
\end{lem}

\begin{proof} 
For a covering $\UC$ of $Y$ consisting of mutually disjoint 
subsets that are closed
and open -- let's call it a cut-and-paste covering -- the $R$-submodule
  $\FC_\UC(Y;R)\subset\FC(Y;R)$ is defined by 
\[\FC_\UC(Y;R)=\{f\in\FC(X;R);~\text{$f_{\vert U}$ constant for all
  $U\in\UC$}\}.\]
Note that a cut-and-paste covering has finitely many elements because
$Y$ is compact. 
Obviously, $\FC_\UC(Y;R)$ is isomorphic to $\bigoplus_{\UC} R$, 
in particular it is a flat $R$-module. 
Further $\FC(Y;R)$ is the union 
$\bigcup_{\UC}\FC_\UC(Y;R)$ where $\UC$ runs through all 
cut-and-paste coverings. Note that this is a directed union; the
coverings form a directed set with respect to the "being
finer"-relation. Since directed colimits and, in particular, directed
unions of flat modules, are again flat~\cite[prop.~(4.4) on
p.~123]{lam(1999)}, we are done. 
\end{proof}

From that we easily deduce the following crucial fact, which will 
be used throughout this paper. 

\begin{lem}\label{flatness-of-crossed-product-ring}
Let $R$ be a ring $R$, and be $B$ a compact subset of 
a topological space $Y$. Then $\FC(Y;R)\rtimes\Gamma\chi_B$ is a 
flat $R\Gamma$-module. 
\end{lem}

\begin{proof} 
Just note that
$M\otimes_{R\Gamma}\FC(Y;R)\rtimes\Gamma\chi_B\cong M\otimes_R\FC(B;R)$ 
and use lemma~\ref{flatness-of-ring-of-functions}. 
\end{proof}

\begin{rem}\label{remark on Morita equivalence} 
We recall some basics about Morita 
equivalences (cf.~\cite[(18.30) on p.~490]{lam(1999)}). 
Let $R$ be a ring and be $p\in R$ an idempotent, i.e.~$p^2=p$. 
Then $p$ is called \textbf{full} if $RpR=R$ holds. In this case, 
the rings $pRp$ and $R$ are \textbf{Morita equivalent} 
and tensoring with the 
bimodules $pR$ resp.~$Rp$ provides equivalences of their respective 
module categories. Furthermore, these equivalences are exact and 
preserve projective modules. 
\end{rem}

\begin{lem}\label{restriction-to-compact-subset-preserves-projectives} 
Let $Y$ be a topological space on which $\Gamma$ acts continuously, 
and be $A\subset Y$ a compact-open subset which meets every 
$\Gamma$-orbit, i.e.~$\Gamma A=Y$. Let $P$ be a projective 
left $R(Y\rtimes\Gamma)$-module. Then $\chi_AP$ is a projective module 
over $\chi_AR(Y\rtimes\Gamma)\chi_A=R(A\rtimes\Gamma)$. 
Furthermore, if $Y$ is compact, then $\chi_A$ is a full idempotent 
in $R(Y\rtimes\Gamma)$, hence the rings $R(Y\rtimes\Gamma)$ and 
$R(A\rtimes\Gamma)$ are Morita equivalent. 
\end{lem}

\begin{proof} 
Recall that we identify $\FC(Y;R)\rtimes\Gamma$ with 
$R(Y\rtimes\Gamma)$. 
For the first assertion, 
it suffices to treat the case $P=\FC(Y;R)\rtimes\Gamma$. 
By assumption we have $Y=\bigcup_{\gamma\in\Gamma}\gamma A$. 
Since the open-closed sets form a set algebra, there 
are open-closed $A_\gamma\subset A$, $\gamma\in\Gamma$, such that 
$Y=\bigcup_{\gamma\in\Gamma} \gamma A_\gamma$ is a disjoint union. 
Note that 
\[\chi_A\FC(Y;R)\rtimes\Gamma\chi_{\gamma A_\gamma}\cong\chi_A\FC(Y;R)\rtimes
\Gamma\chi_{A_\gamma}=\chi_A\FC(Y;R)\rtimes\Gamma\chi_A\chi_{A_\gamma}\]
is a left projective $\chi_A\FC(Y;R)\rtimes\Gamma\chi_A$-module. We
claim that the 
map 
\[\omega:
\chi_A\FC(Y;R)\rtimes\Gamma\longrightarrow\bigoplus_{\gamma\in\Gamma}\chi_A\FC(Y;R)
\rtimes\Gamma\chi_{\gamma A_\gamma}, x\mapsto (x\chi_{\gamma
  A_\gamma})_{\gamma\in\Gamma}\]
is an isomorphism of $\chi_A\FC(Y;R)\rtimes\Gamma\chi_A$-modules. The 
only thing to check is that $\omega$ is well defined: We have to show 
that for $x\in\chi_A\FC(Y;R)\rtimes\Gamma$ there is a finite subset 
$F(x)\subset\Gamma$ such that 
$x\chi_{\gamma A_\gamma}=0$ for $\gamma\not\in F(x)$. 
Consider $x=\sum_{i=1}^nf_i\gamma_i$ with $f_i=f_i\chi_A$. 
Define $F(x)=\{\gamma\in\Gamma;~\exists i=1,\ldots,n:~A\cap
\gamma_i\gamma A_\gamma\ne\emptyset\}$. The set $F(x)$ is 
finite since for any fixed $i\in\{1,\ldots,n\}$ the compact 
set $A$ is covered by the family of disjoint open sets
$\gamma_i\gamma A_\gamma$, $\gamma\in\Gamma$. So we obtain 
for $\gamma\not\in F(x)$ 
\[\sum_{i=1}^nf_i\gamma_i\chi_{\gamma
  A_\gamma}=\sum_{i=1}^nf_i\chi_A\gamma_i\chi_{\gamma
  A_\gamma}=\sum_{i=1}^nf_i\chi_{A\cap\gamma_i\gamma A}\gamma_i=0.\]
If, additionally, $Y$ is compact, there is a finite 
subset $F\subset\Gamma$ with $Y=\bigcup_{\gamma\in F}\gamma A_\gamma$,
  $A_\gamma\subset A$ as above. But then we get $1=\chi_Y=\sum_{\gamma\in F}
  \gamma\chi_A\chi_{A_\gamma}\gamma^{-1}$, so $\chi_A$ is full. 
\end{proof}

\section{The Induction Map in Group (Co)homology} 
\label{sec:the transfer map in group homology} 

Throughout this section \textbf{let $\phi:\Lambda\rightarrow\Gamma$ be a
uniform embedding} (definition~\ref{definition-of-uniform-embedding}). 
By theorem~\ref{dynamic-criterion} there is a topological coupling 
$X$ on which $\Gamma$ and 
$\Lambda$ act in a commuting way, and closed-open fundamental domains 
$X_\Gamma,X_\Lambda$ with \textbf{$X_\Gamma$ being compact}. 
If $\phi$ is a quasi-isometry, then such an $X$ exists with both 
$X_\Lambda,X_\Gamma$ being compact. 
Further, there are compact-open subsets 
$A\subset X_\Gamma, B\subset X_\Lambda$ and a groupoid isomorphism 
$\Phi:A\rtimes\Lambda\overset{\cong}{\rightarrow} B\rtimes\Gamma$ 
(lemma~\ref{isomorphism-of-transformation-groupoids}) inducing a ring 
isomorphism
$\Phi :R(B\rtimes\Gamma)\overset{\cong}{\rightarrow}R(A\rtimes\Lambda)$ 
(lemma~\ref{isomorphism-between-groupoid-rings}), where $R$ is an
arbitrary commutative 
ring. \textbf{We fix a choice of $X, X_\Gamma,X_\Lambda, A, B$ and 
$\Phi$ for the rest of this section}. 
Recall the identification
$\chi_AR(X_\Gamma\rtimes\Lambda)\chi_A=R(A\rtimes\Lambda)$
(remark~\ref{restricted-groupoid-ring}). Further it will be 
crucial in the sequel because of its algebraic consequences 
(see~\ref{remark on Morita equivalence}) 
that the characteristic function $\chi_A$ is 
a full idempotent in $R(X_\Gamma\rtimes\Lambda)$
by~\ref{restriction-to-compact-subset-preserves-projectives}. If 
$X_\Lambda$ is also compact, i.e.~in the case of a quasi-isometry, 
$\chi_B$ is also a full idempotent in $R(X_\Lambda\rtimes\Gamma)$. 
\\[\smallskipamount]
Next we define an induction map from the group homology
resp.~cohomology of $\Lambda$
to that of $\Gamma$. This induction map \textit{depends on the 
topological coupling and on $\Phi$},  
but we suppress that in the notation to simplify it. The reader should
be warned that because of 
these choices it is \textit{not} functorial with 
respect to compositions of uniform embeddings. 
First we 
define induction functors on the level of modules. 
If $f:R\rightarrow S$ is a ring homomorphism and 
$M$ an $S$-module, we denote by $\res_f M$ the $R$-module $M$ induced
by $f$. If $f$ is the inclusion of 
a subring $R\subset S$ we also write $\res_R M$. 
Consider the functors 
\begin{gather}\label{eq:homological-induction-on-the-level-of-modules} 
\iu:\{\text{$R\Lambda$-modules}\}\longrightarrow
\{\text{$R\Gamma$-modules}\}\\
 M\longmapsto\iu(M)=\res_{R\Gamma}\Bigl(R(X_\Lambda\rtimes\Gamma)\chi_B
\otimes_{R(B\rtimes\Gamma)}\res_\Phi\bigl(\chi_AR(X_\Gamma\rtimes\Lambda)
\otimes_{R\Lambda} M\bigl)\Bigl)\notag
\end{gather}
and 
\begin{gather}\label{eq:cohomological-induction-on-the-level-of-modules} 
\io:\{\text{$R\Lambda$-modules}\}\longrightarrow
\{\text{$R\Gamma$-modules}\}\\
 M\longmapsto\io(M)=\res_{R\Gamma}\Bigl(\hom_{R(B\rtimes\Gamma)}\bigl(
\chi_BR(X_\Lambda\rtimes\Gamma),\res_\Phi\bigl(\chi_AR(X_\Gamma\rtimes\Lambda)\otimes_{R\Lambda}M\bigr)\Bigr)\notag  
\end{gather}
Explicitly, the (left) $R\Gamma$-module structure on $\io(M)$ is 
given by $(\gamma f)(x)=f(x\gamma)$,
$x\in\chi_BR(X_\Lambda\rtimes\Gamma)$. 

\begin{lem} 
Assume that $\phi:\Lambda\rightarrow\Gamma$ is a quasi-isometry, 
in particular $X_\Lambda$ is compact. 
Then $\io(R)$ and $\iu(R)$ are 
both isomorphic to $\FC(X_\Lambda)=\FC(X_\Lambda;R)$ as $R\Gamma$-modules. 
\end{lem}

\begin{proof}
For $\iu(R)$ this follows from 
\begin{align*}
\iu(R) &=\res_{R\Gamma}\Bigl(R(X_\Lambda\rtimes\Gamma)\chi_B
\otimes_{R(B\rtimes\Gamma)}\res_\Phi\bigl(\chi_AR(X_\Gamma\rtimes\Lambda)
\otimes_{R\Lambda} R\bigl)\Bigl)\\
&= \res_{R\Gamma}\Bigl(R(X_\Lambda\rtimes\Gamma)\chi_B
\otimes_{R(B\rtimes\Gamma)}\res_\Phi\FC(A)\Bigl)\\
&= \res_{R\Gamma}\Bigl(R(X_\Lambda\rtimes\Gamma)\chi_B
\otimes_{R(B\rtimes\Gamma)}\FC(B)\Bigl)\\
&= \res_{R\Gamma}\Bigl(R(X_\Lambda\rtimes\Gamma)
\otimes_{R(X_\Lambda\rtimes\Gamma)}\FC(X_\Lambda)\Bigl)\\
&= \res_{R\Gamma}\FC(X_\Lambda).
\end{align*}
The fourth equality comes from the fact that 
$\chi_B$ is a full idempotent in $R(X_\Lambda\rtimes\Gamma)$ if 
$X_\Lambda$ is compact. 
Similarly, one concludes $\io(R)\cong\FC(X_\Lambda)$. 
\end{proof}

\begin{thm}[Induction for Uniform
Embeddings]\label{theorem-induction-in-homology}  
Let $M$ be an $R\Lambda$-module. 
Then there are homomorphisms, natural in $M$, 
\[I_n: \xymatrix@1{\homologie_n(\Lambda,M)\ar[r] &\homologie_n(\Gamma,\iu(M))},~n\ge 0,\]
between the group homology of $\Lambda$ and $\Gamma$ with 
coefficients in $M$ resp.~$\iu(M)$. 
We say $I_n$ is the \textbf{induction homomorphism}  
(with respect to the topological coupling X and the groupoid map
$\Phi$ above ). The map $I_n$ factorizes as 
\[\xymatrix@1{\homologie_n(\Lambda,M)\ar[r]&
\homologie_n(\Lambda,\res_{R\Lambda}R(X_\Gamma\rtimes\Lambda)\otimes_{R\Lambda}
M)\ar[r]^-{\cong} & \homologie_n(\Gamma,\iu(M))},\]
where the first map is induced by the map of the coefficients 
$j:M\rightarrow R(X_\Gamma\rtimes\Lambda)\otimes_{R\Lambda} M, m\mapsto
1\otimes m$, and the second map is 
an isomorphism. Similarly, 
there are natural homomorphisms 
\[I^n: \xymatrix@1{\homologie^n(\Lambda,M)\ar[r]& \homologie^n(\Gamma,\io(M))},~n\ge 0,\]
factorizing as $\xymatrix@1{\homologie^n(\Lambda,M)\ar[r]&
\homologie^n\bigl(\Lambda,\res_{R\Lambda}R(X_\Gamma\rtimes\Lambda)\otimes_{R\Lambda}M
\bigr)\ar[r]^-{\cong} &\homologie^n\bigl(\Gamma,\iu(M)\bigr)},$
where again the first map is induced by $j$ 
and the second map is an isomorphism. Furthermore, if $\phi$ is an 
quasi-isometry, in particular $X_\Lambda$ is compact, then 
the induction homomorphisms for $M=R$ 
\begin{gather*}
\homologie_\ast(\Lambda,R)\rightarrow
\homologie_\ast\bigl(\Gamma,\FC(X_\Lambda)\bigr)\\
\homologie^\ast(\Lambda,R)\rightarrow\homologie^\ast\bigl(\Gamma,\FC(X_\Lambda
)\bigr)
\end{gather*}
are compatible with cup- and cap-products. 
\end{thm}

\begin{proof}
The proof occupies the rest of this section. We prefer to be slightly 
redundant in separating the argument for the compatibility 
of the product structures from the rest, thereby improving  
the readability.\\[\smallskipamount] 
Let $P_\ast$ be a right projective $R\Lambda$-resolution of the 
trivial module $R$, and be $Q_\ast$ a right projective 
$R\Gamma$-resolution of
$R$. By lemma~\ref{flatness-of-crossed-product-ring} 
$R(X_\Gamma\rtimes\Lambda)$ is flat over $R\Lambda$ and 
$R(X_\Lambda\rtimes\Gamma)\chi_B$ is flat over $R\Gamma$. 
In particular,
$P_\ast\otimes_{R\Lambda}R(X_\Gamma\rtimes\Lambda)$ and
$Q_\ast\otimes_{R\Gamma}R(X_\Lambda\rtimes\Gamma)\chi_B$ are 
$R(X_\Gamma\rtimes\Lambda)$- resp.~$R(B\rtimes\Gamma)$-resolutions of 
$\FC(X_\Gamma;R)$ resp.~$\FC(B;R)$. Both are projective resolutions; 
for the first this obvious, for the latter this is due 
to lemma~\ref{restriction-to-compact-subset-preserves-projectives}. 
In particular, if $M,N$ are left resp.~right
$R(X_\Gamma\rtimes\Lambda)$-modules, then we have a canonical 
isomorphism $N\otimes_{R(X_\Gamma\rtimes\Lambda)} M\cong
N\chi_A\otimes_{R(A\rtimes\Lambda)}\chi_AM$ since $\chi_A$ is a full
idempotent. Now the composition of the 
following chain maps induces $I_n$, $n\ge 0$, in homology. 
\[\xymatrix{
P_\ast\otimes_{R\Lambda} M\ar[d]^-{1}\\
P_\ast\otimes_{R\Lambda}\bigl(R(X_\Gamma\myrtimes\Lambda)\otimes_{R\Lambda}
M\bigr)\ar[d]^-{2}\\
\bigl(P_\ast\otimes_{R\Lambda}R(X_\Gamma\myrtimes\Lambda)\bigr)\otimes_{R(X_
  \Gamma 
  \rtimes\Lambda)}\bigl(R(X_\Gamma\myrtimes\Lambda)\otimes_{R\Lambda}
M\bigr)\ar[d]^-{3}\\
\bigl(P_\ast\otimes_{R\Lambda}R(X_\Gamma\myrtimes\Lambda)\chi_A\bigr)\otimes_{R(A
  \rtimes\Lambda)}\bigl(\chi_AR(X_\Gamma\myrtimes\Lambda)\otimes_{R\Lambda}
M\bigr)\ar[d]^-{4}\\
\res_\Phi\Bigl(\bigl(P_\ast\otimes_{R\Lambda}R(X_\Gamma\myrtimes\Lambda)\chi_A
\bigr)\Bigr)\otimes_{R(B\rtimes\Gamma)}\res_\Phi\Bigl(\chi_AR(X_\Gamma\myrtimes\Lambda)
\otimes_{R\Lambda}M\Bigr)\ar[d]^-{5}\\
Q_\ast\otimes_{R\Gamma}\underbrace{R(X_\Lambda\myrtimes\Gamma)\chi_B\otimes_{R(B\rtimes\Gamma)} 
\res_\Phi\Bigl(\chi_AR(X_\Gamma\myrtimes\Lambda)
\otimes_{R\Lambda}M\Bigr)}_{=I(M)} 
}\]
The first map is coming from the inclusion $M\rightarrow
R(X_\Gamma\rtimes\Lambda)\otimes_{R\Lambda} M$, $m\mapsto 1\otimes
m$. The second map and the third map are the canonical
identifications, and the fourth map is the obvious isomorphism coming 
from the ring isomorphism $\Phi:R(B\rtimes\Gamma)\overset{\cong}{\rightarrow}
R(A\rtimes\Lambda)$. Note that
$\res_\Phi\bigl(P_\ast\otimes_{R\Lambda}R(X_\Gamma\rtimes\Lambda)\chi_A\bigr)$
and $Q_\ast\otimes_{R\Gamma}R(X_\Lambda\rtimes\Gamma)\chi_B$ are 
projective resolutions of the same module $\FC(B;R)$. Now 
the fifth map is the homotopy equivalence we get 
from the fundamental theorem of homological algebra, which is unique 
up to chain homotopy.  \\[\smallskipamount]
Now let us turn to the cohomological case. Let $P_\ast$,
$Q_\ast$ be left projective $R\Lambda$- resp.~$R\Gamma$-resolutions of
$R$. We only describe the dual diagram: 
%{\small
\[\xymatrix{\hom_{R\Lambda}(P_\ast,M)\ar[d]^1\\
\hom_{R\Lambda}\bigl(P_\ast,R(X_\Gamma\myrtimes\Lambda)\otimes_{R\Lambda}M\bigr)\ar[d]^2\\ 
\hom_{R(X_\Gamma\rtimes\Lambda)}\bigl(R(X_\Gamma\myrtimes\Lambda)\otimes_{R 
  \Lambda}P_\ast,R(X_\Gamma\myrtimes\Lambda)\otimes_{R\Lambda}M\bigr)
\ar[d]^3\\ 
 \hom_{R(A\rtimes\Lambda)}\bigl(\chi_AR(X_\Gamma\myrtimes\Lambda)\otimes_{R
   \Lambda}P_\ast,\chi_AR(X_\Gamma\myrtimes\Lambda)\otimes_{R\Lambda}M\bigr)
 \ar[d]^4\\
\hom_{R(B\rtimes\Gamma)}\!\Bigl(\!\res_\Phi\bigl(\chi_AR(X_\Gamma\myrtimes\Lambda)\!
 \otimes_{R\Lambda}\!P_\ast\bigr),\res_\Phi\bigl(\chi_AR(X_\Gamma\myrtimes\Lambda)\!  
 \otimes_{R\Lambda}\!M\bigr)\Bigr)\ar[d]^5\\
 \hom_{R(B\rtimes\Gamma)}\Bigl(\chi_BR(X_\Lambda\myrtimes\Gamma)
 \otimes_{R\Gamma}Q_\ast,\res_\Phi\bigl(\chi_AR(X_\Gamma\myrtimes\Lambda)
 \otimes_{R\Lambda}M\bigr)\Bigr)
}\]
%}

The five maps correspond to the maps in the 
homological diagram. Now the last term 
is \textit{canonically} isomorphic 
to $\hom_{R\Gamma}(Q_\ast, \io(M))$ by applying 
two canonical isomorphisms of pairs of adjoint 
functors, of the pair
$\chi_BR(X_\Lambda\rtimes\Gamma)\otimes_{R(X_\Lambda\rtimes\Gamma)}\_$ 
and $\hom_{R(B\rtimes\Gamma)}(\chi_BR(X_\Lambda\rtimes\Gamma),\_)$ 
and of the pair
$R(X_\Lambda\rtimes\Gamma)\otimes_{R\Gamma}\_$ and
$\res_{R\Gamma}(\_)$: 
\begin{multline*}
\hom_{R(B\rtimes\Gamma)}\bigl(\chi_BR(X_\Lambda\myrtimes\Gamma)
 \otimes_{R\Gamma}Q_\ast,\res_\Phi\bigl(\chi_AR(X_\Gamma\myrtimes\Lambda)
 \otimes_{R\Lambda}M\bigr)\bigr)\\
\begin{aligned}
\cong &
\hom_{R(X_\Lambda\rtimes\Gamma)}\Bigl(R(X_\Lambda\myrtimes\Gamma)\otimes_{R\Gamma}Q_\ast,\hom_{R(B\rtimes\Gamma)}\bigl(\chi_BR(X_\Lambda\myrtimes\Gamma),
\res_\Phi\bigl(\chi_AR(X_\Gamma\myrtimes\Lambda)\otimes_{R\Lambda}M\bigr)\Bigr)\\
\cong &  \hom_{R\Gamma}\Bigl(Q_\ast,\res_{R\Gamma}
 \underbrace{\hom_{R(B\rtimes\Gamma)}\bigl(\chi_BR(X_\Lambda\myrtimes\Gamma),\res_\Phi\bigl(\chi_AR(X_\Gamma\myrtimes\Lambda)\otimes_{R\Lambda}M\bigr)}_{=\io(M)}\Bigr)   
\end{aligned}
\end{multline*}
Now we turn to the proof of the compatibility with cup- and
cap-products in the case of a quasi-isometry, i.e. $X_\Lambda$ is 
compact. In that case $\chi_B$ is a full idempotent which enables 
us to write the induction map in the following way. 
\[\xymatrix{
H^\ast(\Lambda,R)=\ext_{R\Lambda}^\ast(R,R)\ar[d]_{(\ref{eq:functoriality
of ext})}\\
\ext_{R(X_\Gamma\rtimes\Lambda)}^\ast\bigl(\FC(X_\Gamma),\FC(X_\Gamma)\bigr)
\\
\ext_{R(A\rtimes\Lambda)}^\ast\bigl(\FC(A),\FC(A)\bigr)\ar[u]^{~(\ref{eq:isomorphism between ext-groups coming from a Morita equivalence})}_\cong\ar[d]^\cong\\
\ext_{R(B\rtimes\Gamma)}^\ast\bigl(\FC(B),\FC(B)\bigr)\ar[d]_{(\ref{eq:isomorphism between ext-groups coming from a Morita equivalence})}^\cong\\
H^\ast\bigl(\Gamma,\FC(X_\Lambda)\bigr)=\ext_{R(X_\Lambda\rtimes\Gamma)}^\ast\bigl(\FC(X_\Lambda),\FC(X_\Lambda)\bigr)}
\]
At this point we need the product structures on 
Tor- and Ext-Groups explained in the
appendix. By lemma~\ref{compatibility-of-composition-and-cup-product1} the
cup and the composition product on $H^\ast(\Lambda,R)$ coincide. The 
homomorphisms~(\ref{eq:functoriality of ext}) 
and~(\ref{eq:isomorphism between ext-groups coming from a 
Morita equivalence}) defined in the appendix are compatible with
the composition product. The isomorphism in the middle comes from
the ring isomorphism $\Phi:R(B\rtimes\Gamma)\rightarrow
R(A\rtimes\Lambda)$ and the fact that $\res_\Phi\FC(A)\cong\FC(B)$, 
so it also respects the composition product. 
The composition and the cup product on the last 
term $H^\ast(\Gamma,\FC(X_\Lambda))$ coincide by
lemma~\ref{two-product-structures-on-group-cohomology-coincide}. 
This shows the compatibility with
respect to the cup product. 
The induction in homology can be 
expressed analogously in terms of Tor and the maps~(\ref{eq:functoriality
  of tor}) and~(\ref{eq:isomorphism between tor-groups coming from a
    Morita equivalence}). Since these maps are compatible with respect 
to the evaluation product and since the evaluation and cap product
coincide on the first and last term of the composition by the same lemma, 
the proof is now completed. 
\end{proof}

\section{Quasi-Isometry and (Co)homological Dimension} 
\label{sec:quasi-isometry and homological dimension} 

In this section we prove 
theorem~\ref{statement-theorem estimate with homological dimension}. 
Recall that the \textbf{homological dimension} $\hd_R(\Gamma)$ 
of a group $\Gamma$ 
over a ring $R$ is defined as 
\[\hd_R(\Gamma)=\sup\{n;~\exists\text{$R\Gamma$-module $M$ with
  $\homologie_n(\Gamma,M)\ne 0$}\}\in\NN\cup\{\infty\}.\]
In the same way one defines the cohomological dimension
$cd_R(\Gamma)$. It is a basic fact in group homology
(\cite[lemma 4.1.10]{Weibel(1994)},~\cite[p.~185]{Brown(1994a)})
that $\hd_R(M)$ is the minimal number $n$ such that there 
is a resolution of the trivial $R\Gamma$-module $R$ 
\[0\leftarrow R\leftarrow F_0\leftarrow F_1\leftarrow\cdots\leftarrow
F_n\leftarrow 0\]
by \textit{flat} 
$R\Gamma$-modules $F_i$. Analogously, the cohomological dimension 
$\cd_R(\Gamma)$ is the minimal $n$ such that there is a 
\textit{projective} $R\Gamma$-resolution of $R$ of length
$n$. \\[\smallskipamount] 
We 
prove theorem~\ref{statement-theorem estimate with homological dimension} 
only for the homological dimension; it is an easy matter of dualizing 
to prove the cohomological statement. 
Let us turn to statement (i). 
Suppose that $\Lambda$ uniformly embeds into $\Gamma$, and assume 
$n=\hd_R(\Lambda)<\infty$. Choose an $R\Lambda$-module $M$ with 
$\homologie_n(\Lambda,M)\ne 0$ and a flat $R\Lambda$-resolution
$F_i,~0\le i\le n$, of length $n$ as above. 
The claim would follow from the injectivity of 
the induction map $I_n$ in theorem~\ref{theorem-induction-in-homology}. 
Write $\FC(Y)=\FC(Y;R)$. 
So it suffices to show that for any compact space $Y$ with a 
continuous $\Lambda$-action the homomorphism 
\begin{equation}\label{eq:first map injective}
\homologie_n(\Lambda,M)\longrightarrow \homologie_n(\Lambda, \FC(Y)\rtimes\Lambda\otimes_{R\Lambda} M), 
\end{equation}
which is induced by 
$M\rightarrow\FC(Y;R)\rtimes\Lambda\otimes_{R\Lambda}M$, $m\mapsto 1\otimes m$,
 is injective. But because of 
\begin{align*} 
\homologie_n\bigl(\Lambda,\FC(Y)\myrtimes\Lambda\otimes_{R\Lambda}\!M\bigr)&=\ker\Bigl(F_n\otimes_{R\Lambda}\FC(Y)\myrtimes\Lambda\otimes_{R\Lambda}\! M\rightarrow F_{n-1}\otimes_{R\Lambda}\FC(Y)\myrtimes\Lambda\otimes_{R\Lambda}\! M\Bigr)\\
\homologie_n\bigl(\Lambda,M\bigr)&=\ker\Bigl(F_n\otimes_{R\Lambda}\! M\rightarrow F_{n-1}\otimes_{R\Lambda}\! M\Bigr)
\end{align*}
we are reduced to show that for a flat $R\Lambda$-module $F$ the map 
\[\sigma_F:F\otimes_{R\Lambda} M\longrightarrow F\otimes_{R\Lambda}\FC(Y)\myrtimes\Lambda\otimes_{R\Lambda} M,~x\otimes m\mapsto x\otimes 1\otimes m\] 
is injective. First note that  
$\sigma_{R\Lambda}:M\rightarrow R\Lambda\otimes_{R\Lambda}\FC(Y)\rtimes\Lambda\otimes_{R\Lambda} M=\FC(Y)\otimes_R M$ is injective since 
the $R$-linear map $R\rightarrow\FC(Y)$ (inclusion of constant functions) 
is split by the evaluation $\FC(Y)\rightarrow R$, $f\mapsto f(y_0)$ at 
some base point $y_0$. Hence $\sigma_F$ is injective for any free  
module $F$. By a theorem of Lazard and Govorov~\cite[theorem (4.34) on
p.~134]{lam(1999)} any flat module is the directed colimit 
of free modules and taking directed colimits is an exact
functor~\cite[theorem 2.6.15]{Weibel(1994)}, 
hence $\sigma_F$ is injective for any flat module. 
\smallskip\\
Next we prove statement (ii) 
of theorem~\ref{statement-theorem estimate with homological dimension}. 
Suppose that $\Lambda$ is amenable and that $R$ contains $\QQ$. 
Then the injectivity of~(\ref{eq:first map injective}) 
is obtained as follows. Because of amenability, we can equip $Y$ 
with a $\Lambda$-invariant probability measure.  
By composing that with a $\QQ$-linear map $\RR\rightarrow\QQ$ that 
maps $1$ to $1$, we obtain a signed finitely additive, 
$\QQ$-valued, $\Lambda$-invariant probability measure $\mu$ on $Y$. 
Since a function in $f\in\FC(Y)$ takes only finitely many values, 
we obtain a well defined integration $\int_Yfd\mu$. 
This little trick is taken from~\cite[proof of theorem 1.5]{shalom(2003)}. 
Further, the 
inclusion $R\Lambda\hookrightarrow\FC(Y)\rtimes\Lambda$ is 
split by the $R\Lambda$-\textit{bimodule} map 
\[\FC(Y)\rtimes\Lambda\longrightarrow R\Lambda,~\sum
f_i\lambda_i\mapsto\sum\left(\int_Yf_id\mu\right)\lambda_i.\]
Hence we get a $R\Lambda$-linear left inverse of the map 
$M\rightarrow\FC(Y)\rtimes\Lambda\otimes_{R\Lambda}M$, and 
so~(\ref{eq:first map injective}) is injective. 
This finishes the proof 
of theorem~\ref{statement-theorem estimate with homological dimension}.

\section{Quasi-Isometry and the Cohomology Ring of a Nilpotent
  Group}\label{sec:quasi-isometry and the cohomology ring of a
  nilpotent group} 

\subsection{A module structure on the reduced cohomology}\label{subsec:reduced cohomology} 
We briefly recall the definition of the reduced cohomology 
$\redhomologie^n\!(\Gamma,V)$ of a discrete group $\Gamma$ with coefficients 
in a unitary or orthonormal 
representation $V$ of $\Gamma$. It is defined in terms 
of the \textbf{standard homogeneous resolution}. Consider 
the chain complex 
\begin{equation}\label{eq:standard homogeneous resolution}
C^n(\Gamma,V)=\{\omega:\Gamma^{n+1}\rightarrow
V;~\omega(\gamma\gamma_0,\ldots,\gamma\gamma_n)=\gamma\omega(\gamma_0,\ldots,\gamma_n)\}
\end{equation}
equipped with the standard homogeneous differential 
$d^n:C^n(\Gamma,V)\rightarrow C^{n+1}(\Gamma,V)$
\[(d^n\omega)(\gamma_0,\ldots,\gamma_{n+1})=\sum_{i=0}^{n+1}(-1)^i\omega(\gamma_0,\ldots,\hat{\gamma_i},\ldots,\gamma_{n+1}).\]
Then $C^n(\Gamma,V)$ carries the topology of pointwise convergence. 
The space of $n$-cocycles $\ker d^n$ is 
closed with respect to this topology, but the space of boundaries 
$\im d^{n-1}$
need not be a closed subspace. The \textbf{reduced cohomology} 
$\redhomologie^n\!(\Gamma,V)$ is defined by kernel modulo the closure of the 
image: 
\[\redhomologie^n\!(\Gamma,V)=\frac{\ker\bigl(d^n:C^n(\Gamma,V)\rightarrow
  C^{n+1}(\Gamma,V)\bigr)}{\closure\Bigl(\im\bigl(d^{-1}:C^{n-1}(\Gamma,V)\rightarrow
  C^n(\Gamma, V)\bigr)\Bigr)}\]
There is an obvious surjection
  $\homologie^n(\Gamma,V)\rightarrow\redhomologie^n\!(\Gamma,V)$. \\[\smallskipamount]
Let $Y$ be a compact space equipped with a $\Gamma$-invariant finite 
Borel measure $\mu$, and denote by $L^2(Y;\mu)$ the $\mu$-square-integrable 
real functions on $Y$. We write $\FC(Y)$ for $\FC(Y;\RR)$. 
Since $\FC(Y)$ is a ring
(multiplication of functions) with a $\Gamma$-equivariant
multiplication, the cohomology $H^\ast(Y,\FC(Y))$ is a graded ring by its 
cup-product (see also the appendix). 
In terms of the standard homogeneous resolution the product of 
two cocycles $f:\Gamma^{n+1}\rightarrow\FC(Y)$, 
$g:\Gamma^{m+1}\rightarrow\FC(Y)$ is explicitly given by 
\begin{align}\label{eq:multiplication-formula}
[f]\cup [g]=[\Gamma^{m+n+1}\rightarrow
\FC(Y),(\gamma_0,\ldots,\gamma_{n+m})\mapsto
f(\gamma_0,\ldots,\gamma_n)g(\gamma_n,\ldots,\gamma_{n+m})]. 
\end{align}
We now exhibit a module structure on the reduced 
cohomology $\redhomologie^\ast(Y;L^2(Y;\mu))$ that turns out to be crucial 
for the proof of theorem~\ref{thm-statement: betti numbers of
nilpotent groups}. 
The multiplication of functions 
defines a $\Gamma$-equivariant map 
$\omega:\FC(Y)\otimes_\RR L^2(Y;\mu)\rightarrow
L^2(Y;\mu)$ and hence a graded left 
$\homologie^\ast(\Gamma,\FC(Y))$-module 
structure on $\homologie^\ast(\Gamma,L^2(Y;\mu))$. 
In terms of the standard homogeneous resolution~(\ref{eq:standard 
homogeneous resolution}), this module structure is explicitly given
by the same formula as above: 
\begin{align}\label{eq:module-formula}
[f]\cup [g]=[\Gamma^{m+n+1}\rightarrow
L^2(Y;\mu),(\gamma_0,\ldots,\gamma_{n+m})\mapsto
f(\gamma_0,\ldots,\gamma_n)g(\gamma_n,\ldots,\gamma_{n+m})]
\end{align}
where $f:\Gamma^{n+1}\rightarrow\FC(Y)$, $g:\Gamma^{m+1}\rightarrow
L^2(Y;\mu)$ are cocycles (see~\cite[(11.11) on
p.~65]{Guichardet(1980)}. 
From that we see its continuity for fixed $f$, and so it 
descends to a left \textbf{$\homologie^\ast(\Gamma,\FC(Y))$-module
  structure on the reduced cohomology}
$\redhomologie^\ast\!(\Gamma,L^2(Y;\mu))$. 
Similarly, we get a right $\homologie^\ast(\Gamma,\FC(Y))$-module
structure  
on $\redhomologie^\ast\!(\Gamma,L^2(Y;\mu))$. 

\subsection{Proof of theorem~\ref{thm-statement: betti numbers of
    nilpotent groups}}\label{subsec:proof of theorem on
  cohomology ring of nilpotent groups} 

We now prove 
theorem~\ref{thm-statement: betti numbers of nilpotent groups} 
which says that the real cohomology ring of a finitely generated 
nilpotent group is a quasi-isometry invariant. 
The proof follows from the following theorem, which is a 
"multiplicative" generalization of Shalom's~\cite[theorem
4.1.1]{shalom(2003)} and its method of proof.

\begin{thm}\label{multiplicative-generalization-of-shaloms-betti-result} 
Let $\Gamma$ be a finitely generated amenable group with the 
following property: Any unitary $\Gamma$-representation $\pi$ 
with $\redhomologie^n(\Gamma,\pi)\ne 0$ contains the trivial 
representation. Let $\Lambda$ be a finitely generated 
group which is quasi-isometric to $\Gamma$, and assume that the 
Betti numbers 
$b_n(\Lambda), b_n(\Gamma)$ are finite for all $n\ge 0$. Then there is a 
multiplicative injective homomorphism 
\[\xymatrix@1{\homologie^\ast(\Lambda,\RR)\,\ar@{^(->}[r] &\homologie^\ast(\Gamma,\RR)}\]
between the real cohomology rings of $\Lambda$ and $\Gamma$. 
\end{thm}

First recall the well known fact 
that the classifying spaces of finitely generated torsionfree nilpotent 
groups are finite, and a finitely generated nilpotent group 
is virtually finitely generated 
torsionfree nilpotent, so the hypothesis on the 
Betti numbers is satisfied. 
The fact that finitely generated nilpotent 
groups satisfy the representation-theoretic 
assumption of the theorem is proved  
in~\cite[theorem 4.1.3]{shalom(2003)}. It is deduced from 
a theorem of 
Blanc which says that the continuous 
cohomology of a connected nilpotent Lie group with 
coefficients in an irreducible, non-trivial 
unitary representation 
always vanishes. So this is the only place where 
the nilpotence hypothesis of theorem~\ref{thm-statement: betti numbers
of nilpotent groups} plays an essential role. 
Let $\Lambda,\Gamma$ be finitely generated, nilpotent. 
Due to symmetry (and like in Shalom's article), 
theorem~\ref{multiplicative-generalization-of-shaloms-betti-result}
yields that the Betti numbers of finitely 
generated nilpotent groups are the same. So the injective map 
$H^\ast(\Lambda;\RR)\rightarrow H^\ast(\Gamma;\RR)$ in the 
preceding theorem must be an isomorphism. This completes 
the proof of~\ref{thm-statement: betti numbers of nilpotent groups} 
once we have shown the theorem above.

\begin{proof}[Proof of
theorem~\ref{multiplicative-generalization-of-shaloms-betti-result} ]  
By theorem~\ref{dynamic-criterion} there 
is a topological coupling $X$ with compact fundamental 
domains $X_\Lambda, X_\Gamma$, and by 
theorem~\ref{ergodic ume coupling out of a qi} we can equip
$X$ with a $\Gamma\times\Lambda$-invariant non-trivial 
\textit{ergodic} measure. 
By replacing $\Gamma$ by a product $\Gamma\times F$ with a finite
group
(remark~\ref{shaloms-topological-coupling-version-with-inclusion}) 
we can assume that $X_\Gamma\subset X_\Lambda$. Note that this does
not affect the assumption of the theorem. 
We write $\FC(X_\Gamma)$ for $\FC(X_\Gamma;\RR)$ and $L^2(X_\Gamma)$
for the real $\mu$-square-integrable functions. In the sequel we 
use Shalom's explicit formula~\cite[section
3.2,~(12)]{shalom(2003)} 
for the induction rather than the 
abstract setup of section~\ref{sec:the transfer map in group
homology}, since we work with the reduced cohomology. 
In terms of the standard homogeneous resolution Shalom's induction map
$\homologie^\ast(\Lambda,\RR)\rightarrow
\homologie^\ast(\Gamma,L^2(X_\Lambda))$ is given by 
sending a cycle 
$w:\Lambda^{n+1}\rightarrow\RR$ 
to $\omega':\Gamma^{n+1}\rightarrow L^2(X_\Lambda)$ 
defined by 
\begin{align}\label{eq:induction-formula}
\omega'(\gamma_0,\ldots,\gamma_n)(x)=\omega(\alpha(\gamma_0,x),\ldots,\alpha(\gamma_n,x)).
\end{align} 
Being continuous it descends to a map 
\[I^\ast_{red}: \redhomologie^\ast(\Lambda,\RR)\rightarrow\redhomologie(\Gamma,L^2(X_\Lambda)).\]
By~\cite[theorem
3.2.1]{shalom(2003)} $I^\ast_{red}$ is injective. 
Further note that 
$\homologie^\ast(\Lambda,\RR)=\redhomologie^\ast(\Lambda,\RR)$ and 
$\homologie^\ast(\Lambda,\RR)=\redhomologie^\ast(\Lambda,\RR)$ hold 
because, by hypothesis, the homology groups are finite-dimensional. 
Since $\{x;\alpha(\gamma,x)=\lambda\}$ is open and closed for 
fixed $\gamma,\lambda$, the functions
$\omega'(\gamma_0,\ldots,\gamma_n)$ actually lie in $\FC(X_\Lambda)$ (compare
with the proof of
lemma~\ref{isomorphism-between-groupoid-rings}). Hence 
$I^\ast_{red}$ factorizes as 
\[\xymatrix@1{\homologie^\ast(\Lambda,\RR)\ar[r]^{I^\ast}&\homologie^\ast(\Gamma,\FC(X_\Lambda))\ar[r]^j&
\redhomologie^\ast(\Gamma,L^2(X_\Lambda)),}\]
where $j^\ast$ is the composition of the inclusion 
of coefficients and the canonical surjection 
$\homologie^\ast(\Gamma,L^2(X_\Lambda))\rightarrow
  \redhomologie^\ast\!(\Gamma,L^2(X_\Lambda))$. 
At this point we want to remark that 
it can be shown that this definition of
$I^\ast$ is the same as the one in 
section~\ref{sec:the transfer map in group homology}, but this not
relevant here since we stick to Shalom's definition. 
Now let us collect all these data in the following diagram:  
\[\xymatrix{
 & \redhomologie^\ast\!(\Gamma,L^2(X_\Lambda))\ar@<4ex>@/^2pc/[dd]^{p_2^\ast}\\
 \homologie^\ast(\Lambda,\RR)\ar@{^(->}[ru]^{I^\ast_{red}}\ar@{^(->}[r]^{I^\ast}\ar[rd]_{I^\ast\circ
  p_1^\ast}&
\homologie^\ast(\Gamma,\FC(X_\Lambda))\ar@/^/[d]^{p_1^\ast}\ar[u]^{j^\ast}\\
 & \homologie^\ast(\Gamma,\RR)\ar@{^(->}[u]^{i^\ast}
}\] 
The map $i^\ast$ is the map induced by the inclusion of coefficients. 
The maps $p_1^\ast,p_2^\ast$ 
are the maps induced by the ($\Gamma$-equivariant!) 
integration maps
$L^2(X_\Lambda)\rightarrow\RR,f\mapsto\int_{X_\Lambda}fd\mu$ 
resp.~$\FC(X_\Lambda)\rightarrow\RR,f\mapsto\int_{X_\Lambda}fd\mu$.\\[\smallskipamount]   
Our aim is to show that the composition $p_1^\ast\circ I^\ast$ 
is multiplicative and injective. 
We have a direct sum $L^2(X_\Lambda)=\RR\oplus L^2_0(X_\Lambda)$
 of $\Gamma$-representations, where $\RR$ are the constant functions and
 $L^2_0(X_\Lambda)$ are the functions with mean value zero. Similarly, 
one obtains a direct sum decomposition 
$L^2(X_\Lambda;\CC)=\CC\oplus L^2_0(X_\Lambda;\CC)$ in the complex
 case and a $\RR\Gamma$-module 
decomposition $\FC(X_\Lambda)=\RR\oplus\FC_0(X_\Lambda)$. Hence we get 
direct sum decompositions
\begin{align*}
\redhomologie^\ast\!(\Gamma,L^2(X_\Lambda))&=\redhomologie^\ast\!(\Gamma,\RR)\oplus\redhomologie^\ast\!(\Gamma,L^2_0(X_\Lambda)),\\
\homologie^\ast(\Gamma,\FC(X_\Lambda))&=\homologie^\ast(\Gamma,\RR)\oplus
 \homologie^\ast(\Gamma,\FC_0(X_\Lambda))
\end{align*}
which are respected by $j^\ast$. 
For an element $x\in \homologie^m(\Gamma,\FC(X_\Lambda))$ we use the 
notation 
$x=x_0+x_1,~x_0\in \homologie^m(\Gamma,\FC_0(X_\Lambda)),
~x_1\in \homologie^n(\Gamma,\RR)$ for the sum decomposition. 
Since the $\Gamma$-action on $X_\Lambda$ is 
ergodic, the $\Gamma$-representation $L^2_0(X_\Lambda;\CC)$ does not contain 
the trivial representation. So our assumption yields 
that $\redhomologie^n\!(\Gamma,L^2_0(X_\Lambda;\CC))=0$. Since 
$L^2_0(X_\Lambda;\CC)=L^2_0(X_\Lambda)\oplus L^2_0(X_\Lambda)$ as real 
$\Gamma$-representations, it is
$\redhomologie^n\!(\Gamma,L^2_0(X_\Lambda))=0$. In particular,   
$p_2^\ast$ is an isomorphism. Since 
$p_1^\ast\circ I^\ast=p_2^\ast\circ I_{red}^\ast$ and 
$I_{red}^\ast$ is injective, we obtain the injectivity of 
$p_1^\ast\circ I^\ast$.  \\[\smallskipamount] 
From the formulas~(\ref{eq:multiplication-formula})
and~(\ref{eq:induction-formula}) one can see that  
$I^\ast$ is multiplicative, so 
we are reduced to show that 
the map $p_1^\ast$ is multiplicative in our special 
situation (in general, it is not since it is induced by
integration!). The map $j$ is obviously a module
homomorphism with respect to both the left and right
$\homologie^\ast(\Gamma,\FC(X_\Lambda))$-module structure on
$\homologie^\ast(\Gamma,\FC(X_\Lambda))$ resp.~on
$\redhomologie^\ast(\Gamma,L^2(X_\Lambda))$. For the module 
structure see the preceding section. 
Since $p_2^\ast$ is an isomorphism, it follows that 
$\ker p_1^\ast=\ker j^\ast$, and hence $\ker p_1^\ast$ is a two-sided 
ideal of the cohomology ring $\homologie^\ast(\Gamma,\FC(X_\Lambda))$. 
Further note $\ker p_1^\ast=H^\ast(\Gamma,\FC_0(X_\Lambda))$. 
Now let $x\in\homologie^m(\Gamma,\FC(X_\Lambda))$,
$x'\in\homologie^n(\Gamma,\FC(X_\Lambda))$, and consider 
the sum decompositions $x=x_0+x_1, x'=x_0'+x_1'$ as above. 
It is clear that $p_1^{m+n}(x_1\cup x_1')=p_1^m(x_1)\cup p_1^n(x_1')$ 
since $i^\ast$, induced by the ring homomorphism
$\RR\rightarrow\FC(X_\Lambda)$, is a multiplicative map. 
Now one computes: 
\begin{align*} 
p_1^{m+n}(x\cup x') &=p_1^{m+n}\bigl((x_0+x_1)\cup (x_0'+x_1')\bigr)\\
&= p_1^{m+n}(x_0\cup x_0')+p_1^{m+n}(x_0\cup x_1')+p_1^{m+n}(x_1\cup
x_0')+p_1^{m+n}(x_1\cup x_1')\\
&= p_1^{m+n}(x_1\cup x_1')\\
&= p_1^{m}(x_1)\cup p_1^n(x_1')\\
&= p_1^m(x)\cup p_1^n(x').
\end{align*}
\end{proof}

\subsection{Remarks}\label{subsec:examples}

The $\QQ$-Malcev completions of quasi-isometric nilpotent 
groups are in general not isomorphic. In fact, they are only isomorphic 
if the groups are commensurable. This fact 
is also reflected in cohomology: We 
give examples of quasi-isometric nilpotent 
groups whose rational cohomology rings are not isomorphic 
as $\QQ$-algebras. So 
theorem~\ref{thm-statement: betti numbers of nilpotent groups} 
does not hold for rational coefficients. Consider 
the $2n+1$-dimensional Heisenberg Lie algebra $H(n,K)$ over a field 
$K$. The is the Lie algebra of matrices with entries in $K$ of the form 
\[\left(\begin{array}{ccc} 0 & x & z\\
                      0 & 0 & y\\
                      0 & 0 & 0
        \end{array}\right)\]
where $x$, $y$ are an $n$-dimensional row resp.~column vector and $z\in K$. 
Now consider the $4n+2$-dimensional, 
\textit{rational} Lie algebras $\lieg=H(n,\QQ(\sqrt{2}))$ and 
$\lieh=H(n,\QQ)\oplus H(n,\QQ)$. Over $\QQ$, $\lieg$ and $\lieh$ are 
not isomorphic since the $\QQ$-codimension of the centralizers 
of elements in $\lieg$ is either $0$ or $\ge 2$, whereas in $\lieh$ 
there is an element whose centralizer has $\QQ$-codimension $1$. 
However, $\lieg\otimes_\QQ\RR$ and $\lieh\otimes_\QQ\RR$ are 
easily seen to be isomorphic as real Lie algebras. This example 
is discussed in~\cite[remark 2.15]{raghunathan(1972)}. 
Hence $\lieg$, $\lieh$ are two different rational structures of the 
same real Lie algebra, thereby defining two quasi-isometric (but not 
commensurable) nilpotent Lie groups $\Lambda$ and $\Gamma$. By virtue
of Nomizu's theorem~\cite{nomizu(1954)} $\homologie^\ast(\Lambda,\QQ)$ and 
$\homologie^\ast(\lieg,\QQ)$ resp.~$\homologie^\ast(\Gamma,\QQ)$ and 
$\homologie^\ast(\lieh,\QQ)$ are isomorphic as algebras. 
Now assume that $n\ge 2$. 
Due to~\cite[example 2.3]{carlson-toledo(1995)} 
the Chevalley-Eilenberg complexes of 
$H(n, K)$, hence also of 
$H(n,K)\oplus H(n,K)$, as differential graded algebras 
are $1$-formal, in the sense of rational homotopy theory. 
Note 
that $1$-formality is independent of the ground field. 
The reason behind this is that these algebras are quadratically 
presented. In particular, the rational Lie algebras $\lieg$, $\lieh$ can be 
reconstructed from their rational cohomology algebras (compare
with~\cite[prop.~2.1]{carlson-toledo(1995)}. Therefore the latter cannot 
be isomorphic. At this point we remark that the class of quadratically 
presented Lie algebras is quite restricted within the class of all 
nilpotent Lie algebras. See~\cite{carlson-toledo(1995)} for more 
information. \\[\medskipamount]
Next we discuss what information we can get if the representation 
theoretic assumption of
theorem~\ref{multiplicative-generalization-of-shaloms-betti-result} is 
not satisfied in every dimension but only for the first reduced 
cohomology. Explicitly, we assume only that any unitary
$\Gamma$-representation $\pi$ with $\redhomologie^1(\Gamma,\pi)\ne 0$ 
contains the trivial representation. This property of $\Gamma$ was 
called \textit{property $H_T$} in~\cite[definition of
p.~7]{shalom(2003)}. Otherwise we retain the same assumptions on 
$\Lambda$ and $\Gamma$ as 
in~\ref{multiplicative-generalization-of-shaloms-betti-result}. 
By the same arguments we get homomorphisms 
$\alpha^n:H^n(\Lambda,\RR)\rightarrow\homologie^n(\Gamma,\RR)$, where
$\alpha^1$ is injective. Further, if $x\in\homologie^1(\Lambda,\RR)$, 
$y\in\homologie^n(\Lambda,\RR)$, $n\ge 0$, then 
\[\alpha^{n+1}(x\cup
y)=\alpha^1(x)\cup\alpha^n(y)\] 
holds true. To show that, 
we must verify that the integration map $p_1^\ast$ in the proof 
satisfies 
\begin{equation}\label{eq:mult}
p_1^{n+1}(z\cup w)=p_1^1(z)\cup p_1^n(w)
\end{equation}
for $z\in\homologie^1(\Gamma,\FC(X_\Lambda)), w\in\homologie^n(\Gamma,
\FC(X_\Lambda))$. Decompose $z=z_0+z_1$ as in the proof. Then one 
can use the module structure on the reduced cohomology to conclude 
$p_1^{n+1}(z_0\cup w)=0$. On the other hand, 
$p_1^{n+1}(z_1\cup w)=p_1^1(z_1)\cup p_1^n(w)$ since 
$p_1$ is obviously a $H^\ast(\Gamma,\RR)$-module homomorphism. 
Hence~(\ref{eq:mult}) follows. 

\section{Quasi-Isometry Invariance of Novikov-Shubin Invariants} 
\label{sec:quasi-isometry invariance of novikov-shubin invariants}

\subsection{Review of the definition of Novikov-Shubin
  invariants}\label{subsec:review of the definition of novikov-shubin
  invariants} 
In this section we review the usual definition of 
Novikov-Shubin invariants and then the 
\textit{homological viewpoint} developed
in~\cite{Lueck-Reich-Schick(1999)}. 
The latter is strongly  
inspired by L\"uck's algebraic approach to $L^2$-Betti
numbers. \\[\smallskipamount]
Let $\Gamma$ be a discrete group. The Hilbert space with basis
$\Gamma$ is denoted by $l^2(\Gamma)$. 
The bounded operators on 
$l^2(\Gamma)$ which are equivariant with respect to the left 
$\Gamma$-action on $l^2(\Gamma)$ form a von Neumann algebra
$\NC(\Gamma)$, called 
the \textbf{group von Neumann algebra}. Equivalently, $\NC(\Gamma)$ 
can be defined as the weak closure of $\CC\Gamma$ whose elements 
act on $l^2(\Gamma)$ by right multiplication. 
A von Neumann algebra $\AC$ is 
called \textbf{finite} if $\AC$ comes equipped with a \textbf{finite 
trace}, i.e.~with a finite, normal, faithful trace
$\trace_\AC:\AC\rightarrow\CC$. The examples of finite von Neumann
algebras we will consider are the following. 
\begin{enumerate}[$\centerdot$]
\item The group von Neumann algebra $\NC(\Gamma)$ with its standard
  trace $\trace_{\NC(\Gamma)}:\NC(\Gamma)\rightarrow\CC$ given by 
  $\trace_{\NC(\Gamma)}(T)=\scalar{1_\Gamma}{T(1_\Gamma)}_{l^2(\Gamma)}$. In
  particular, the trace of an element in $\CC\Gamma$ is the 
  coefficient of the unit element. 
\item Let $X$ be a standard Borel space equipped with a probability 
Borel measure $\mu$. Then $\lxb=L^\infty(X;\mu)$ 
is a finite von Neumann algebra 
with the trace $\trace_\lxb(f)=\int_Xfd\mu$. 
\item Assume $X$ is additionally equipped with a $\mu$-preserving 
$\Gamma$-action. Then there is the \textbf{von Neumann crossed
  product} $\lxb\neurtimes\Gamma$ which contains the
\textit{algebraic} crossed product $\lxb\rtimes\Gamma$ as a weakly 
dense subalgebra. Further, $\lxb\neurtimes\Gamma$ has a finite trace 
whose restriction to $\lxb\rtimes\Gamma$ is given by 
$\trace(\sum_{\gamma}f_\gamma\gamma)=\int_X f_1d\mu$. 
\item If $\AC$ is a finite von Neumann algebra with 
trace $\trace_\AC$, then the $n$-dimensional 
square matrices $M_n(\AC)$ over $\AC$ form 
again a finite von Neumann algebra with the trace 
\[\trace_{M_n(\AC)}\bigl((T_{ij})_{1\le i,j\le
  n}\bigr)=\frac{1}{n}\sum_{i=1}^n\trace_\AC\bigl(T_{ii}\bigr).\]
\end{enumerate}

Any finite von Neumann algebra $\AC$ is contained in an  
algebra of unbounded operators which 
we explain next. A not necessarily bounded operator $a$ on $\BC(H)$ is 
\textbf{affiliated} to $\AC\subset\BC(H)$ if $ba\subset ab$ for all operators 
$b\in\AC'$ holds. Here $\AC'$ is the commutant of $\AC$, and
  $ba\subset ab$ means that restricted to the possibly smaller domain
  of $ba$ the two operators coincide. 
\begin{defi}\label{defintion of affiliated operators}
Let $\UC(\AC)$ denote the set of all closed densely defined operators 
affiliated to $\AC$. If $\AC$ happens to be $\NC(\Gamma)$ 
we write $\UC(\Gamma)$ instead of $\UC(\NC(\Gamma))$. 
\end{defi} 

Murray and von Neumann~\cite{murray-neumann(1936)}    
showed that these unbounded operators form a $\ast$-algebra 
containing $\AC$, when addition, multiplication and involution are
defined as the closures of the naive addition, multiplication and 
involution. More information about 
$\UC(\AC)$ can be found in~\cite[chapter 8]{Lueck(2002)}
and~\cite{Reich(2001)}. From the latter reference we take also 
the following result. 

\begin{thm}\label{flatness of affiliated operators}
The inclusion $\AC\subset\UC(\AC)$ is a flat 
ring extension. 
\end{thm}

\begin{defi}\label{definition of novikov-shubin invariant of an
    operator} 
Let $\AC$ be a finite von Neumann algebra with trace $\trace_\AC$, and 
for a normal operator $T\in\AC$ 
we denote by $f(T)$ the operator obtained from spectral calculus with 
respect to the function $f$. The 
\textbf{Novikov-Shubin invariant $\alpha_\AC(T)\in
[0,\infty]\cup\infty^+$ of a positive operator $T\in\AC$} is defined as 
\[\alpha_\AC(T)=\begin{cases}\displaystyle\liminf_{\lambda\rightarrow
    0^+}\frac{\ln\bigl (\trace_\AC (\chi_{(0,\lambda]}(T))\bigr)}{\ln
    (\lambda)}\in [0,\infty]&\text{ if $\trace_\AC
    (\chi_{(0,\lambda]}(T))>0$ for $\lambda>0$,}\\
    \infty^+&\text{ otherwise.}
  \end{cases}\]
Here $\chi_{(0,\lambda]}$ is the characteristic function of the 
interval $(0,\lambda]$. 
\end{defi}

Now we explain the usual definition of the Novikov-Shubin
invariant of a discrete group $\Gamma$ having a finite type model of 
its classifying space $\textup{E}\Gamma$. The chain complex 
of such a model is a 
free $\ZZ\Gamma$-resolution of $\ZZ$ 
\[0\leftarrow\ZZ\leftarrow F_0\overset{d_1}{\longleftarrow} 
F_1\overset{d_2}{\longleftarrow}F_2\overset{d_3}{\longleftarrow}\cdots,\]
where for each $i\ge 0$ 
the module $F_i$ is a finitely generated free $\ZZ\Gamma$-module, let
us say of rank $r(i)$. Now by
tensoring this complex with $l^2(\Gamma)$ we get a Hilbert complex 
$l^2(\Gamma)\otimes_{\ZZ\Gamma}F_\ast$ whose operators are bounded 
$\Gamma$-equivariant operators $\bar{d}_\ast$. 
Then, by choosing a basis for each $F_i$, we can consider 
$\bar{d}_i^*\bar{d}_i$ as an element in $M_{r(i)}(\NC(\Gamma))$. 
Now the following makes sense and is indeed 
well defined (see~\cite[chapter 2]{Lueck(2002)} for proofs and much more
information). 

\begin{defi}[classical definition]\label{definition
    novikov-shubin invariants classical}
The $i$-th Novikov-Shubin invariant $\alpha_i(\Gamma)$, $i\ge 1$, 
of $\Gamma$ is defined as
$\alpha_{M_{r(i)}(\NC(\Gamma))}(\bar{d}_i^*\bar{d}_i)$. 
\end{defi}

For the homological viewpoint one defines 
the notion of capacity for any $\NC(\Gamma)$-module in the first
place, and then the Novikov-Shubin invariants of 
$\Gamma$ are essentially the reciprocals of the capacities of the 
group homology of $\Gamma$ with coefficients in the $\CC\Gamma$-module 
$\NC(\Gamma)$. To a module $M$ 
over a finite von Neumann algebra $\AC$ there are 
two numerical invariants attached, the \textbf{dimension} 
$\dim_{\AC}(M)\in [0,\infty]$ in the sense of 
L\"uck~\cite{Lueck(1998a),Lueck(1998b)} and 
the \textbf{capacity} 
$c_{\AC}(M)\in [0,\infty]\cup\{0^-\}$~\cite{Lueck-Reich-Schick(1999)}. 
We will not elaborate upon neither the dimension nor the question of 
$c_{\AC}(M)$ being well defined; the reader is 
referred to consult~\cite{Lueck(2002)} for the first
and~\cite{Lueck-Reich-Schick(1999)} for the latter. 
The definition of $c_\AC(M)$ proceeds in three steps.

\begin{enumerate}[\sc 1. Step:]
\item Let $M$ be a finitely 
 presented $\AC$-module $M$ with $\dim_\AC(M)=0$. 
By~\cite[lemma 3.4, definition 3.11]{Lueck(1997a)} 
there is a short exact sequence 
\begin{center}
$\xymatrix@1{\AC^n\,\ar@{^(->}[r]^{r_A}&\AC^n\ar@{->>}[r]& M}$,
\end{center}
where $r_A$ is right multiplication with a positive matrix 
$A\in M_n(\AC)$. Define 
\begin{center}
$\displaystyle c_\AC(M)=\frac{1}{\alpha_{M_n(\AC)}(A)}\in
[0,\infty]\cup\{0^-\}$, 
\end{center}
where we agree on the following rules 
\begin{equation*}
1/\infty=0,~1/\infty^+=0^-\text{ and }0^-<0.
\end{equation*} 
\item An $\AC$-module $M$ is called \textbf{measurable} if it is the quotient 
of a finitely presented $\AC$-module $N$ with $\dim_\AC (N)=0$. 
For such a module $M$ we define 
$c^{'}_\AC(M)$ as 
\begin{center}
$c^{'}_\AC(M)=\inf\{c_\AC(N);~\text{$N$ finitely presented, $\dim_\AC(N)=0$,
  $N$ surjects onto $M$}\}$.
\end{center}
\item For an arbitrary $\AC$-module $M$ we set 
\begin{center}
$c^{''}_\AC(M)=\sup\{c_\AC(N);~\text{$N$ measurable submodule of
  $M$}\}$. 
\end{center}
In~\cite[proposition 2.4]{Lueck-Reich-Schick(1999)} it is shown that 
for a finitely presented $\AC$-module $M$ with $\dim_\AC(M)=0$ and a 
measurable $\AC$-module $N$ one has $c_\AC(M)=c^{'}_\AC(M)$ and 
$c^{'}_\AC(N)=c^{''}_\AC(N)$. Therefore we do not need to distinguish
between $c_\AC, c^{'}_\AC, c^{''}_\AC$. We call $c_\AC(M)$ the
\textbf{capacity} of $M$. 
\end{enumerate}

\begin{defi}\label{definition of cofinal-measurable} 
An $\AC$-module $M$ is \textbf{cofinal-measurable} if each finitely 
generated submodule is measurable. This is equivalent to 
$\UC(\AC)\otimes_\AC M=0$~\cite[note 4.5]{Reich(2001)}. 
\end{defi}

\begin{defi}[homological approach] Let $\Gamma$ be an arbitrary 
discrete group. Then we define the \textbf{$i$-th capacity $c_i(\Gamma)$ of 
$\Gamma$} as 
\[c_i(\Gamma)=c_{\NC(\Gamma)}\bigl(\homologie_i(\Gamma,\NC(\Gamma))\bigr), i\ge 
0.\]
The \textbf{$i$-th Novikov-Shubin invariant $\alpha_i(\Gamma)$, $i\ge 1$, 
of $\Gamma$} is defined as the inverse of $c_{i-1}(\Gamma)$ with the 
usual rules $1/0=\infty$, $1/0^-=\infty^+$. 
\end{defi}

Of course, this can be shown to be compatible with the classical 
definition whenever it applies. What is the point about 
introducing a name for the reciprocal of the Novikov-Shubin 
invariant? In some sense, $c_i(\Gamma)$ measures the \textit{size} of 
$\homologie_i(\Gamma,\NC(\Gamma))$ if 
the $i$-th $L^2$-Betti number vanishes. In particular, it equals $0^-$ if 
$\homologie_i(\Gamma,\NC(\Gamma))=0$. For this reason, it is more convenient 
in the homological setting to work with the inverse. 
As already indicated in the introduction, 
the definition of the capacity is not 
well-behaved for arbitrary $\Gamma$. When we restrict to 
amenable groups, here is the subclass (see theorem~\ref{description of
  the class CM} for its description) where $c_i(\Gamma)$  
still behaves well:

\begin{defi}\label{definition of the class CM} 
Let $\CM$ be the class of all finitely generated 
  amenable groups $\Gamma$ such that $\homologie_n(\Gamma,\NC(\Gamma))$ is
  cofinal-measurable for $n\ge 1$. Equivalently, 
we can say that $\CM$ consists precisely of all 
finitely generated amenable groups $\Gamma$ with 
$\homologie_n(\Gamma,\UC(\Gamma))=0$ for $n\ge 1$. 
\end{defi}

\subsection{Some technical properties of the capacity}\label{subsec:some
  properties of the capacity}

\begin{thm}\label{capacity-induction}
 Let $\AC\subset\BC$ be a trace-preserving 
inclusion of finite von Neumann algebras, and let $M$ be an 
$\AC$-module. Then the following holds. 
 \begin{enumerate}[(i)]
\item The ring extensions $\AC\subset\BC$ and
  $\UC(\AC)\subset\UC(\BC)$ are faithfully
  flat. \label{capacity-induction i)}
\item If $M$ is a cofinal-measurable $\AC$-module, 
then $\BC\otimes_\AC M$ is cofinal-measurable and 
$c_\AC(M)=c_\BC(\BC\otimes_\AC M)$ holds. \label{capacity-induction ii)}
 \item If $M$ is a finitely presented 
 $\AC$-module, then $\BC\otimes_\AC M$ is finitely presented and 
$c_\AC(M)=c_\BC(\BC\otimes_\AC M)$ holds. \label{capacity-induction iii)}
 \end{enumerate}
 \end{thm}

 \begin{proof}\hfill\\
(i) In~\cite[6.29]{Lueck(2002)} the first statement in (i)  
is proved for inclusions of  
group von Neumann algebras induced by inclusions of groups. 
We give a (shorter) proof of the general case. 
We show that $\BC$ is a torsionless 
$\AC$-module, i.e.~for every $b\in\BC$, $b\ne 0$, there exists 
an $\AC$-linear map $f:\BC\rightarrow\AC$ with $f(b)\ne 0$. 
Any von Neumann algebra is 
a semihereditary ring~\cite[remark on p.~288]{Lueck(2002)},  
and every torsionless module 
over a semihereditary ring is flat~\cite[theorem 4.67]{lam(1999)}, 
hence this implies that $\AC\subset\BC$ is flat.  
Let $b$ be a non-zero element in $\BC$.  
Since $b$ is the sum of four unitaries~
\cite[theorem 4.1.7 on p.~242]{kadison-ringrose(1997a)}, there is a  
unitary $u\in\BC$ such that $\trace_\BC(b^*u)\ne 0$. The map 
$\AC\rightarrow\BC$,  
$a\mapsto au$ extends to an $\AC$-equivariant isometric embedding 
$i:l^2\AC\rightarrow l^2\BC$ between the GNS-representations of
$\trace_\AC$ resp.~$\trace_\BC$. Taking the orthogonal 
projection onto the image of $i$ yields an $\AC$-equivariant bounded 
split $f:l^2\BC\rightarrow l^2\AC$ of $i$. 
Since 
$\scalar{b}{u}_{l^2\BC}=\trace (b^*u)\ne 0$ we get $f(b)\ne 0$. 
It remains to show that $f(\BC)\subset\AC$. 
Due to a theorem of Dixmier~\cite[theorem 9.9]{Lueck(2002)}, 
this follows from the fact that 
$\AC\rightarrow l^2\AC$, $n\mapsto n\cdot f(m)=f(n\cdot m)$ 
extends to a bounded operator on $l^2\AC$. 
The construction of $f$ for $u=1$ defines an $\AC$-equivariant 
split of the inclusion $\AC\subset\BC$. Hence $\AC\subset\BC$ is 
a faithfully flat ring extension. 
The algebra of affiliated operators $\UC(\AC)$ is a von Neumann 
regular ring~\cite[2.1]{Reich(2001)} (see
also~\cite[theorem 8.22]{Lueck(2002)}). Hence every $\UC(\AC)$-module 
is flat~\cite[theorem 4.21]{lam(1999)}. This implies also 
that every ring extension $\UC(\AC)\subset R$ into a bigger ring 
$R$ is faithfully flat~\cite[theorem 4.74
(4)]{lam(1999)}. \\[\smallskipamount] 
(ii) This is proved in~\cite[lemma 2.12]{Lueck-Reich-Schick(1999)} 
for inclusions of group von Neumann algebras but the proof 
carries over to arbitrary inclusions without
modifications. \\[\smallskipamount]  
(iii) From the exactness of $\AC\subset\BC$ follows immediately that 
$\BC\otimes_\AC M$ is again finitely presented. 
We use the fact that a finitely presented $\AC$-module $M$ splits as 
$M=\boldP M\oplus \boldT M$ 
into a finitely generated projective module 
$\boldP M$ and into a finitely presented 
module $\boldT M$ with $\dim_\AC (\boldT M)=0$~\cite[theorem
0.6]{Lueck(1998a)}. 
The capacity satisfies $c_\AC(M)=c_\AC(\boldT M)$ 
by~\cite[theorem 2.7]{Lueck-Reich-Schick(1999)}. 
Similarly, we obtain 
$c_\BC(\BC\otimes_\AC M)=c_\BC(\BC\otimes_\AC \boldT M)$. 
Since $\boldT M$ is cofinal-measurable, (iii) follows now from (ii). 
\end{proof}

\begin{thm}
 Let $\AC$ be a finite von Neumann algebra and $p\in\AC$ a full 
 projection. Then for every $p\AC p$-module $M$ the equality 
$c_\AC(\AC p\otimes_{p\AC p}M)=c_{p\AC p}(M)$ holds. 
\end{thm}

\begin{proof}
First assume that $M$ is zero-dimensional and finitely presented. 
Then there exists a short exact sequence~\cite[lemma
3.4]{Lueck(1997a)} 
of the form 
\[\xymatrix@1{\bigl (p\AC p\bigr )^n\,\ar@{^(->}[r]^{r_A}&\bigl
(p\AC p\bigr )^n\ar@{->>}[r]& M}.\]
where $r_A$ is right multiplication with a positive matrix 
$A\in M_n(p\AC p)$. Tensoring by $\AC p$ is flat (see
remark~\ref{remark on Morita equivalence}), so we obtain a short exact sequence
 \[\xymatrix@1{\bigl (\AC p\bigr )^n\oplus\bigl (\AC (1-p)\bigr
   )^n\,\ar@{^(->}[r]^-{r_A\oplus id}&\bigl (\AC
   p\bigr )^n\oplus\bigl (\AC\bigl (1-p)\bigr )^n\ar@{->>}[r]& \AC
   p\otimes_{p\AC p}M}.\]
The map in the middle is an $\AC$-equivariant endomorphism of $\AC^n$ 
given by right multiplication with the matrix 
$A+(1-p)\id_{\AC^n}\in M_n(\AC)$. In particular, we have 
that $\AC p\otimes_{p\AC p} M$ is again zero-dimensional since 
the dimension $\dim_\AC$ is additive~\cite[theorem 0.6]{Lueck(1998a)}. 
Now let us show two general facts about spectral calculus 
in a von Neumann algebra $\AC$. 
Let $f:\RR\rightarrow\RR$ be an essentially bounded Borel function 
with $f(0)=0$. \\
\noindent (1) For selfadjoint $a,b\in\AC$ 
satisfying $a\cdot b=b\cdot a=0$ we have the following 
equality of operators obtained by spectral calculus in $\AC$ 
with respect to $f$. 
\[f(a+b;\AC)=f(a;\AC)+f(b;\AC)\]
\noindent (2) For a projection $p\in\AC$ and selfadjoint $a\in p\AC p$ the
equality  
\[f(a;\AC)=pf(a;\AC)p=f(a;p\AC p)\]
holds. 
Since $f(a;\AC)$ is the strong limit of operators of the form
$q(a)$ where $q$ is a polynomial with vanishing constant term, it 
suffices to show these equalities for $f(x)=x^n$, $n\ge 1$. 
This is obvious; note for (2) that $pa=ap=a$ holds. 
Notice that we have 
$\chi_{(0,\lambda]}((1-p)\id_{\AC^n})=0$ for $\lambda<1$, 
because $(1-p)\id_{\AC^n}$ is a projection. Now we compute using (1)
and (2) for $\lambda<1$
\begin{align*}
\trace_{M_n(\AC)} \bigl (\chi_{(0,\lambda]}(A+(1-p)\id_{\AC^n};\AC)\bigr )
  &=\trace_{M_n(\AC)} \bigl (\chi_{(0,\lambda]}(A;\AC)\bigr)\\
  &=\trace_{M_n(\AC)} \bigl (\chi_{(0,\lambda]}(A;p\AC p)\bigr )\\
  &=\trace_{\AC}(p)\cdot\trace_{M_n(p\AC p)} \bigl
  (\chi_{(0,\lambda]}(A;p\AC p)\bigr ).
\end{align*}
Thus $\alpha_{M_n(p\AC p)}(A)=\alpha_{M_n(\AC)} (A+(1-p)\id_{\AC^n})$
holds, hence the assertion for $M$ follows. 
\smallskip\\
Now assume that $M$ is measurable. For $\epsilon>0$ let 
$N$ be a zero-dimensional,
finitely presented $p\AC p$-module surjecting onto $M$ such that $c_{p\AC
  p}(N)\le c_{p\AC p}(M)-\epsilon$. As seen above, 
$\AC p\otimes_{p\AC p} N$ 
is again finitely presented zero-dimensional, and it surjects onto 
$\AC p\otimes_{p\AC p} M$. 
Then we obtain 
\[c_\AC(\AC p\otimes_{p\AC p}M)\le c_\AC (\AC p\otimes_{p\AC
  p}N)=c_{p\AC p} (N)\le c_{p\AC p}(M)-\epsilon,\]
thus $c_\AC(\AC p\otimes_{p\AC p}M)\le c_{p\AC p}(M)$. 
Since the Morita equivalence $\AC p\otimes_{p\AC p}\_$ has the inverse 
$p\AC\otimes_\AC\_$ we can conclude analogously to obtain 
$c_{p\AC p}(M)\le c_\AC(\AC p\otimes_{p\AC p}M)$, thus 
$c_{p\AC p}(M)=c_\AC(\AC p\otimes_{p\AC p}M)$. The argument for an 
arbitrary module $M$ is similar. 
\end{proof}

\begin{cor}\label{capacity-idempotent}
Let $\AC$ be a finite von Neumann algebra and $p\in\AC$ a full 
projection. Then for every $\AC$-module $M$ the equality 
$c_\AC(M)=c_{p\AC p}(pM)$ holds. 
\end{cor}

\subsection{The proof of theorem~\ref{statement-theorem
    qi invariance of novikov-shubin invariants} and the class
    $\CM$}\label{subsec:proof of qi invariance}   
\begin{proof}[Proof of theorem~\ref{statement-theorem
    qi invariance of novikov-shubin invariants}]
Let $\Gamma$ be a group in $\CM$ and $\Lambda$ be quasi-isometric 
to $\Gamma$. According to~\ref{description of
  the class CM} (ii) $\Lambda$ lies then also in $\CM$. 
In particular both groups are amenable. 
According to theorems~\ref{dynamic-criterion} 
and~\ref{ergodic ume coupling out of a qi} 
there exists a 
topological coupling $X$ of 
$\Gamma,\Lambda$ with compact fundamental domains 
$X_\Gamma,X_\Lambda$ and a non-trivial 
$\Gamma\times\Lambda$-invariant Borel measure $\mu$ on $X$ such 
that $X_\Gamma$, $X_\Lambda$ have finite measure. Due to 
lemma~\ref{isomorphism-between-groupoid-rings} there exist 
open and closed subspaces $A\subset X_\Gamma$, $B\subset X_\Lambda$ 
for which there is a 
trace-preserving ring isomorphism between 
$\chi_A\FC(X_\Gamma;\CC)\rtimes\Lambda\chi_A$ and 
$\chi_B\FC(X_\Lambda;\CC)\rtimes\Gamma\chi_B$. Further, 
$\chi_A$ resp.~$\chi_B$ is a full idempotent in its respective 
crossed product ring and hence in any bigger ring according 
to lemma~\ref{restriction-to-compact-subset-preserves-projectives}. 
From the construction of $X_\Gamma$ as a closed-open subset 
of some space of functions between discrete sets it is clear that 
the closed-open sets of $X_\Gamma$ form a basis of the topology. 
Hence $\FC(X_\Gamma;\CC)\subset L^\infty(X_\Gamma)$ is 
weakly dense. The same applies to $X_\Lambda$. 
Since the isomorphism above is trace-preserving, 
it extends to an isomorphism 
\[\chi_AL^\infty(X_\Gamma)\rtimes\Lambda\chi_A\rightarrow\chi_BL^\infty(
    X_\Lambda)\rtimes\Gamma\chi_B.\]
For the same reason this extends to a trace-preserving isomorphism of the 
respective von Neumann algebras and to an isomorphism of their algebras of
    affiliated operators. 
So we come up with 
the following commutative diagram with vertical isomorphisms. 
\begin{equation}\label{eq:diagram-trace-preserving-isomorphisms} 
\xymatrix{L^\infty(A)\ar@{^(->}[r]\ar[d]^\cong &
  \chi_AL^\infty(X_\Gamma)\rtimes\Lambda\chi_A\ar@{^(->}[r]\ar[d]^\cong& 
  \chi_AL^\infty(X_\Gamma)\neurtimes\Lambda\chi_A\ar[d]^\cong\ar@{^(->}[r]&
  \UC\bigl(\chi_AL^\infty(X_\Gamma)\neurtimes\Lambda\chi_A\bigr)\ar[d]^\cong\\ 
L^\infty(B)\ar@{^(->}[r]& 
\chi_BL^\infty(X_\Lambda)\rtimes\Gamma\chi_B\ar@{^(->}[r] &
\chi_BL^\infty(X_\Lambda)\neurtimes\Gamma\chi_B\ar@{^(->}[r] & 
\UC\bigl(\chi_BL^\infty(X_\Lambda)\neurtimes\Gamma\chi_B\bigr)}
\end{equation} 
Hence it suffices to prove that 
\begin{align}\label{eq:essential capacity equation} 
c_n(\Lambda)&=c_{\chi_AL^\infty(X_\Gamma)\neurtimes\Lambda\chi_A}\Bigl(
\tor_n^{\chi_AL^\infty(X_\Gamma)\rtimes\Lambda\chi_A}\bigl(\chi_A
L^\infty(X_\Gamma)\neurtimes\Lambda\chi_A,L^\infty(A)\bigr)\Bigr)\\
c_n(\Gamma)&=c_{\chi_BL^\infty(X_\Lambda)\neurtimes\Gamma\chi_B}\Bigl(
\tor_n^{\chi_BL^\infty(X_\Lambda)\rtimes\Gamma\chi_B}\bigl(\chi_B
L^\infty(X_\Lambda)\neurtimes\Gamma\chi_B,L^\infty(B)\bigr)\Bigr),\notag
\end{align}
which follows by applying (in this
order) theorem~\ref{capacity-induction}~(\ref{capacity-induction
ii)}), 
theorem~\ref{capacity-induction}~(\ref{capacity-induction
  i)}),~(\ref{eq:adjunction-isomorphism-tor}) and 
~(\ref{eq:isomorphism between tor-groups coming from a Morita
  equivalence}) from the appendix and corollary~\ref{capacity-idempotent}:
\begin{align*} 
c_n(\Lambda)=&c_{\NC(\Lambda)}\Bigl(\homologie_n\bigl(\Gamma,\NC(\Lambda)\bigr)
\Bigr) \\
=&c_{\NC(\Lambda)}\Bigl(\tor_n^{\CC\Lambda}\bigl(\NC(\Lambda),\CC\bigr)
\Bigr) \\
=&c_{L^\infty(X_\Gamma)\neurtimes\Lambda}\Bigl(L^\infty(X_\Gamma)\neurtimes\Lambda\otimes_{\NC(\Lambda)}\tor_n^{\CC\Lambda}\bigl(\NC(\Lambda),\CC\bigr)\Bigr)\\
=&c_{L^\infty(X_\Gamma)\neurtimes\Lambda}\Bigl(\tor_n^{\CC\Lambda}\bigl(L^\infty(X_\Gamma)\neurtimes\Lambda,\CC\bigr)\Bigr)\\
=&c_{L^\infty(X_\Gamma)\neurtimes\Lambda}\Bigl(\tor_n^{L^\infty(X_\Gamma)
  \rtimes\Lambda}\bigl(L^\infty(X_\Gamma)\neurtimes\Lambda,L^\infty(X_\Lambda)\bigr)\Bigr)\\
=&c_{L^\infty(X_\Gamma)\neurtimes\Lambda}\Bigr(\tor_n^{\chi_AL^\infty(X_\Gamma)\rtimes\Lambda\chi_A}\bigl(L^\infty(X_\Gamma)\neurtimes\Lambda\chi_A,L^\infty(A)\bigr)\Bigr)\\
=&c_{\chi_AL^\infty(X_\Gamma)\neurtimes\Lambda\chi_A}\Bigl(
\chi_A\tor_n^{\chi_AL^\infty(X_\Gamma)\rtimes\Lambda\chi_A}\bigl(L^\infty(X_\Gamma)\neurtimes\Lambda\chi_A,L^\infty(A)\bigr)\Bigr)\\
=&c_{\chi_AL^\infty(X_\Gamma)\neurtimes\Lambda\chi_A}\Bigl(
\tor_n^{\chi_AL^\infty(X_\Gamma)\rtimes\Lambda\chi_A}\bigl(\chi_A
L^\infty(X_\Gamma)\neurtimes\Lambda\chi_A,L^\infty(A)\bigr)\Bigr).
\end{align*}
Similar for $\Gamma$. This finishes the proof of
theorem~\ref{statement-theorem qi invariance of novikov-shubin
  invariants}. 
\end{proof}

Basically the same method yields also 
the following result for 
not necessarily amenable groups, which 
should be compared to the 
proportionality theorem in~\cite[theorem 3.183]{Lueck(2002)}. 

\begin{thm}\label{statement-theorem proportionality} 
Let $\Gamma$ and $\Lambda$ be finitely generated groups, 
and assume that $\Gamma$ is of type 
$FP_{n+1}$ over $\CC$. If $\Gamma$ and $\Lambda$ are cocompact lattices 
of the same locally compact group, then 
$\alpha_i(\Gamma)=\alpha_i(\Lambda)$ holds for $1\le i\le n+1$. 
\end{thm}

We indicate how to modify the proof above. Details 
are left to the reader. 
By~\cite[corollary 9]{alonso(1994)} $\Lambda$ is also of 
type $FP_{n+1}$ over $\CC$. This assumption is needed for 
the compatibility of the capacity with induction from the 
group von Neumann algebra to a bigger one
(see theorem~\ref{capacity-induction} (iii)) and replaces the hypothesis
"being in $\Ccal$" needed to apply~\ref{capacity-induction} (ii). 
The topological coupling used for the proof is the common locally 
compact group that contains $\Lambda$, $\Gamma$ as cocompact
lattices.

\begin{thm}\label{description of the class CM}\hfill
\begin{enumerate}[(i)]
\item The class $\CM$ contains all finitely generated 
elementary amenable groups that have a
bound on the orders of their finite subgroups and all 
amenable groups of type $FP_\infty$ over $\CC$. Moreover, if a
finitely generated amenable 
group $\Gamma$ has an infinite, finitely generated, 
normal subgroup in $\CM$, then $\Gamma$ lies in
$\CM$. \label{description of the class CM (i)}
\item The class $\CM$ is geometric, i.e.~if a group
  $\Gamma$ is quasi-isometric to 
a group in $\CM$ then $\Gamma$ lies in $\CM$. \label{description of
  the class CM (ii)} 
\end{enumerate}
\end{thm}

\begin{proof}
(\ref{description of the class CM (i)}) 
Let $\Gamma$ be an amenable group of type $FP_\infty$ over $\CC$. 
Since the category of finitely presented modules over a finite von Neumann
algebra is abelian~\cite[theorem 0.2]{Lueck(1997a)}, it follows that 
$\homologie_n(\Gamma;\NC(\Gamma))$ is a finitely presented 
$\NC(\Gamma)$-module. 
 By a theorem of Cheeger and Gromov
 $\dim_{\NC(\Gamma)}(\homologie_n(\Gamma;\NC(\Gamma)))=0$, $n\ge 1$,  
 holds for amenable $\Gamma$ (see also~\cite[corollary 5.13]{Lueck(1998a)}). 
 By~\cite[8.22 (4)]{Lueck(2002)}, $\dim_{\NC(\Gamma)} (M)=0$ and 
 $\UC(\Gamma)\otimes_{\NC(\Gamma)} M=0$ are equivalent for a finitely
 presented  
 module $M$. This yields
 $\UC(\Gamma)\otimes_{\NC(\Gamma)}\homologie_n(\Gamma;\NC(\Gamma))=0$, $n\ge 
 1$. By theorem~\ref{flatness of affiliated operators} we obtain 
$\homologie_n(\Gamma;\UC(\Gamma))=0$, $n\ge 1$, hence $\Gamma$ lies in the
 class $\CM$.  
 For an elementary amenable group with a bound on the orders of 
 finite subgroups $\UC(\Gamma)$ is a flat $\CC G$-module~\cite[theorem
 9.1]{Reich(1999)}. Hence such a $\Gamma$ lies in $\CM$. 
 Now assume that a group $\Gamma$ has an infinite, finitely generated, 
 normal subgroup $\Lambda$ which lies in $\CM$,
 i.e.~$\homologie_n(\Lambda;\UC(\Lambda))=0$, $n\ge 
 1$. Since $\Lambda$ is finitely generated, the module
 $\homologie_0(\Lambda;\NC(\Lambda))$ is finitely 
 presented, and since $\Lambda$ is infinite, its zeroth $L^2$-Betti number
 vanishes. As seen above, this implies $\homologie_0(\Lambda;\UC(\Lambda))=0$. 
 Since $\UC(\Gamma)$ is flat over
 $\UC(\Lambda)$ by~\ref{capacity-induction} we get $\homologie_n(\Lambda;\UC(\Gamma))=0$
 for $n\ge 0$.  
 The Hochschild-Serre spectral 
 sequence now yields $\homologie_n(\Gamma;\UC(\Gamma))=0$, $n\ge 0$, in particular  
 $\Gamma\in\CM$. \\[\smallskipamount]
(\ref{description of the class CM (ii)}) Let be $\Gamma\in\CM$, and 
$\Lambda$ be quasi-isometric to $\Gamma$. 
The argument is completely analogous to
the one for equation~(\ref{eq:essential capacity equation}): It turns
out that the vanishing of $\homologie_n(\Gamma,\UC(\Gamma))$ is equivalent 
to the vanishing of
\[\tor_n^{\chi_BL^\infty(X_\Lambda)\rtimes\Gamma\chi_B}
\Bigl(\underbrace{\UC\bigl(\chi_BL^\infty(X_\Lambda)\neurtimes\Gamma\chi_B\bigr)}_{=\chi_B\UC(L^\infty(X_\Lambda)\neurtimes\Gamma)\chi_B},L^\infty(B)\Bigr),  
\]  
and similar for $\Lambda$. Due to 
diagram~(\ref{eq:diagram-trace-preserving-isomorphisms}), 
$\homologie_n(\Gamma,\UC(\Gamma))$ vanishes if and only 
if $\homologie_n(\Lambda,\UC(\Lambda))$ vanishes. 
 \end{proof}

\begin{rem}
By a theorem of Cheeger and Gromov~\cite{Cheeger(1986)} 
the $L^2$-Betti numbers of infinite amenable groups 
vanish. However, there are amenable groups $\Gamma$, 
e.g.~the lamplighter group, such that
$\homologie_n(\Gamma;\UC(\Gamma))$ does not vanish 
for some $n\ge 1$. See~\cite{linnell-lueck-schick(2002)}. In
particular, the lamplighter group is \textit{not} in $\CM$. We 
wonder whether the property of belonging 
to $\CM$ is relevant to other spectral issues of amenable groups. 
\end{rem}

\appendix

\section{Product Structures on Tor and Ext}\label{appendix}

Let $R$ be a ring, and be $M$, $N$ $R$-modules. 
Then the Ext-groups $\ext^n_R(M,N)$ can be computed 
with a projective $R$-resolution $P_\ast$ of $M$: 
$\ext^n_R(M,N)=\homologie^n(\hom_R(P_\ast, N))$. 
If $M$ is a right $R$-module, we have 
$\tor_n^R(M,N)=\homologie_n(P_\ast\otimes N)$. 
Notice that $\homologie^n(\Gamma,M)=\ext^n_{R\Gamma}(R,M)$ and 
$\homologie_n(\Gamma,M)=\tor^n_{R\Gamma}(R,M)$. 
We describe now a more \textit{symmetric} 
way to define Ext- and Tor-groups 
which allows to introduce product structures.\\[\smallskipamount]
We recollect some standard notations and facts about chain complexes 
(cf.~\cite[chapter I, 0]{Brown(1994a)}). 
Let $C_\ast$, $C_\ast'$ be non-negative chain complexes of left $R$-modules. 
Let $\hom_R(C_\ast, C_\ast')_n$, $n\ge 0$, be the abelian 
group of graded $R$-module homomorphisms of degree $n$ 
from $C_\ast$ to $C_\ast'$, i.e.
\[\hom_R(C_\ast, C_\ast')_n=\prod_{q\ge 0}\hom_R(C_q,C_{q+n}').\]
Then $\hom_R(C_\ast, C_\ast')_n$, $n\ge 0$, becomes a chain 
complex by the differential 
\begin{gather*} 
\partial_{\hom}:\hom_R(C_\ast,
C_\ast')_n\longrightarrow\hom_R(C_\ast,C_\ast')_{n-1}\\
\partial_{\hom}(f)=\partial f-(-1)^nf\partial.
\end{gather*}
The $0$-cycles of this complex are just the chain maps from $C_\ast$
to $C_\ast'$, and the $0$-boundaries are the nullhomotopic chain
maps. 
Hence its zeroth homology is the abelian group of homotopy classes of
chain maps from $C_\ast$ to $C_\ast'$. If $D_\ast$ is a complex of 
right $R$-modules, then $D_\ast\otimes_R C_\ast$ denotes the tensor product 
of complexes defined by 
\[\bigl(D_\ast\otimes_R C_\ast\bigr)_n=\bigoplus_{i+j=n}D_i\otimes_R C_j\]
and equipped with the differential 
$\partial_{\otimes}(x\otimes y)=\partial x\otimes
y+(-1)^{\operatorname{deg} x}x\otimes\partial y$.   

Now let $M$ and $N$ be left $R$-modules, and be 
$M\leftarrow P_\ast$ and $N\leftarrow Q_\ast$ left projective 
$R$-resolutions. Consider $N$ as a chain complex 
concentrated in degree zero. 
Due to~\cite[theorem (8.5) on p.~29]{Brown(1994a)}, 
the chain map $Q_\ast\rightarrow
N$ induces a quasi-isomorphism 
$\hom_R(P_\ast,Q_\ast)_\ast\rightarrow\hom_R(P_\ast, N)_\ast$, i.e.~the
induced map in homology 
\[\xymatrix@1{H_{-n}(P_\ast,Q_\ast)\ar[r]^-{\cong}& 
\homologie_{-n}(\hom_R(P_\ast,N))}=\homologie^n(\hom_R(P_\ast,N))=\ext^n_R(M,N)
\]  
is a (canonical) isomorphism. Similarly, if $M$ is a right $R$-module
and $N$ is left $R$-module with right resp.~left projective
resolutions $P_\ast$, $Q_\ast$, then the chain map $Q_\ast\rightarrow
N$ induces a quasi-isomorphism $P_\ast\otimes_RQ_\ast\rightarrow
P_\ast\otimes N$, where the homology of the latter is
$\tor_\ast^R(M,N)$. 
This is what we understand under the  
\textit{symmetric way of defining Ext and Tor}. \\[\smallskipamount]
Let us recall some \textbf{functoriality properties of the Ext and Tor}.  
Compare with~\cite[p.~72-74]{Guichardet(1980)}. 
Assume $R\subset S$ is a \textit{flat ring
extension}. Let $M, N$ be $R$-modules (whether right or left will be
clear from the context) with projective 
$R$-resolutions $P_\ast$ resp.~$Q_\ast$. 
Then there are natural maps (\textit{induction of scalars}) 
\begin{gather}
\ext^n_R(N,M)\longrightarrow\ext^n_S(S\otimes_R N,S\otimes_R
M)\label{eq:functoriality of ext},\\
\tor_n^R(N,M)\longrightarrow\tor_n^S(N\otimes_RS,S\otimes_RM)\label{eq:functoriality
  of tor}. 
\end{gather}
In the symmetric picture above,~(\ref{eq:functoriality of ext}) is
induced by scalar induction 
\[\hom_R(Q_\ast,P_\ast)_n\longrightarrow
\hom_r(S\otimes_RQ_\ast, S\otimes_R P_\ast)_n,~f\mapsto S\otimes_R
f.\]
Similar for~(\ref{eq:functoriality
  of tor}). 
If $L$ is an $S$-module then we have -- as a result of the adjunction 
isomorphism between induction and restriction -- canonical 
isomorphisms 
\begin{gather}
\ext^n_S(S\otimes_R M,L)\cong\ext^n_R(M,
L)\label{eq:adjunction-isomorphism-ext}\\
\tor_n^S(M\otimes_RS,L)\cong\tor_n^R(M,L)\label{eq:adjunction-isomorphism-tor}
\end{gather}
If $p$ is a full idempotent in the ring $R$, implying that 
$pRp$ and $R$ are Morita equivalent (cf.~remark~\ref{remark on Morita
  equivalence}), then we obtain natural isomorphisms 
\begin{gather}
\ext^n_{pRp}(N,M)\overset{\cong}{\longrightarrow}\ext^n_{R}(Rp\otimes_{pRp}N,Rp
\otimes_{pRp}M)\label{eq:isomorphism between ext-groups coming from a
    Morita equivalence}\\
\tor_n^{pRp}(N,M)\overset{\cong}{\longrightarrow}\tor_n^{R}(N\otimes_{pRp}pR,Rp
\otimes_{pRp}M)\label{eq:isomorphism between tor-groups coming from a
    Morita equivalence}
\end{gather}
In the symmetric picture, the map~(\ref{eq:isomorphism between ext-groups coming from a
    Morita equivalence}) is given by tensoring chain 
maps with $Rp$ 
\[\hom_R(Q_\ast,P_\ast)_n\longrightarrow
    \hom_R(Rp\otimes_{pRp} Q_\ast, Rp\otimes_{pRp} P_\ast)_n,~f\mapsto
    Rp\otimes_{pRp} f.\]
Its inverse is given by tensoring with $pR$ over $R$. Similar for 
Tor. \\[\smallskipamount]
Next we introduce the multiplicative structures. 
Let $M,N$ and $L$ be $R$-modules with projective $R$-resolutions 
$P_\ast,Q_\ast$ and $W_\ast$. The composition of chain maps 
in the symmetric picture induces homomorphisms 
\[\ext^m_R(N,L)\otimes\ext^n_R(M,N)\longrightarrow\ext^{m+n}_R(M,L)\]
which turn $\ext_R^\ast(M,M)$ into a graded ring whose 
product is called \textbf{composition product}. Further, we have a 
product, called the \textbf{evaluation product}, 
\[\ext_R^m(N,L)\otimes\tor^R_n(M, N)\longrightarrow\tor^R_{n-m}(M,L)\]
defined by 
\[\hom_R(Q_\ast,W_\ast)_m\otimes \bigl(P_\ast\otimes_R
Q_\ast\bigr)_n\rightarrow\bigl(P_\ast\otimes Q_\ast)_{m-n}, f\otimes
p\otimes q\rightarrow (-1)^{\deg f\deg p}p\otimes f(q).\]
Obviously, the homomorphisms 
in~(\ref{eq:functoriality of ext}) and~(\ref{eq:isomorphism between
  ext-groups coming from a Morita equivalence}) are
multiplicative. \\[\smallskipamount] 
Now assume $R$ is \textit{commutative}. Let $M, N, L$ be
$R\Gamma$-modules, and let $c:M\otimes_R N\rightarrow L$ 
be a $R\Gamma$-homomorphism, where $M\otimes_R N$ is equipped with the
diagonal $\Gamma$-action.  
Then the \textbf{cup product} $\cup$ 
and \textbf{cap product} $\cap$ (see~\cite[p.~109-117]{Brown(1994a)}) 
are maps 
\begin{gather*}
\cup: \homologie^m(\Gamma,M)\otimes
\homologie^n(\Gamma,N)\longrightarrow
\homologie^{m+n}(\Gamma,M\otimes_R
N)\overset{c}{\rightarrow}\homologie^{m+n}(\Gamma,L)\\
\cap: \homologie^m(\Gamma,M)\otimes
\homologie^n(\Gamma,N)\longrightarrow
\homologie_{m-n}(\Gamma,M\otimes_R
N)\overset{c}{\rightarrow}\homologie_{m-n}(\Gamma,L) 
\end{gather*}
This yields ring and module structures in the cases $N=M=L$ and $N=L$. 
We record from~\cite[theorem (4.6) on p.~115]{Brown(1994a)}. 
\begin{lem}\label{compatibility-of-composition-and-cup-product1} 
The cup product and composition product resp.~the cap product and 
evaluation product 
  on $\ext^\ast_{R\Gamma}(R,R)=\homologie^\ast(\Gamma, R)$,
  $\tor_\ast^{R\Gamma}(R,R)=\homologie_\ast(\Gamma,R)$ coincide. 
\end{lem}

Let $Y$ be a compact topological space on which a group $\Gamma$ acts. 
Write $\FC(Y)=\FC(Y;R)$. 
By lemma~\ref{flatness-of-ring-of-functions} 
$R\Gamma\subset\FC(Y)\rtimes\Gamma$ is a flat 
ring extension. 
On $\homologie^\ast(\Gamma,\FC(Y))$ we can exhibit 
two product structures. The first one is the cup product coming from 
$\FC(Y)\otimes_R\FC(Y)\rightarrow\FC(Y), f\otimes g\mapsto fg$, 
and the second one is the composition product 
coming from the isomorphism
\[\xymatrix@1{\homologie^n(\Gamma,\FC(Y))=\ext^n_{R\Gamma}(R,\FC(Y))\ar[r]^-{\cong}_-{~(\ref{eq:adjunction-isomorphism-ext})}&\ext^n_{\FC(Y)\rtimes\Gamma}
(\FC(Y) ,\FC(Y))}.\]
Similarly, using~(\ref{eq:adjunction-isomorphism-tor}) 
we obtain a cap and an evaluation product on 
$\homologie^\ast(\Gamma,\FC(Y))$ and $\homologie_\ast(\Gamma,\FC(Y))$. 
Then the following fact follows again from~\cite[theorem (4.6) on
p.~115]{Brown(1994a)}. 
\begin{lem}\label{two-product-structures-on-group-cohomology-coincide}
The cup and composition resp.~the cap and evaluation product on 
$\homologie^\ast(\Gamma,\FC(Y))$, 
$\homologie_\ast(\Gamma,\FC(Y))$ coincide. 
\end{lem}

\bibliographystyle{amsalpha}
\bibliography{promotion,dbdef,dbpre,dbpub}

\def\cprime{$'$} \def\polhk#1{\setbox0=\hbox{#1}{\ooalign{\hidewidth
  \lower1.5ex\hbox{`}\hidewidth\crcr\unhbox0}}}
\providecommand{\bysame}{\leavevmode\hbox to3em{\hrulefill}\thinspace}
\providecommand{\MR}{\relax\ifhmode\unskip\space\fi MR }
\providecommand{\MRhref}[2]{%
  \href{http://www.ams.org/mathscinet-getitem?mr=#1}{#2}
}
\providecommand{\href}[2]{#2}
\begin{thebibliography}{L{\"u}c98b}

\bibitem[Alo94]{alonso(1994)}
Juan~M. Alonso, \emph{Finiteness conditions on groups and quasi-isometries}, J.
  Pure Appl. Algebra \textbf{95} (1994), no.~2, 121--129. \MR{95f:20083}

\bibitem[BG96]{bridson(1996)}
M.~R. Bridson and S.~M. Gersten, \emph{The optimal isoperimetric inequality for
  torus bundles over the circle}, Quart. J. Math. Oxford Ser. (2) \textbf{47}
  (1996), no.~185, 1--23. \MR{97c:20047}

\bibitem[Bro94]{Brown(1994a)}
Kenneth~S. Brown, \emph{Cohomology of groups}, Springer-Verlag, New York, 1994,
  Corrected reprint of the 1982 original. \MR{96a:20072}

\bibitem[BS73]{borelserre(1973)}
A.~Borel and J.-P. Serre, \emph{Corners and arithmetic groups}, Comment. Math.
  Helv. \textbf{48} (1973), 436--491. \MR{52 \#8337}

\bibitem[CG86]{Cheeger(1986)}
Jeff Cheeger and Mikhael Gromov, \emph{{$L\sb 2$}-cohomology and group
  cohomology}, Topology \textbf{25} (1986), no.~2, 189--215.

\bibitem[CT95]{carlson-toledo(1995)}
James~A. Carlson and Domingo Toledo, \emph{Quadratic presentations and
  nilpotent k{\"a}hler groups}, J. Geom. Anal. \textbf{5} (1995), no.~3,
  359--377. \MR{1360825 (97c:32038)}

\bibitem[Ger93]{gersten(1993)}
Steve~M. Gersten, \emph{Quasi-isometry invariance of cohomological dimension},
  C. R. Acad. Sci. Paris S\'er. I Math. \textbf{316} (1993), no.~5, 411--416.
  \MR{94b:20042}

\bibitem[Gro93]{gromov(1993)}
M.~Gromov, \emph{Asymptotic invariants of infinite groups}, Geometric group
  theory, Vol.\ 2 (Sussex, 1991), Cambridge Univ. Press, Cambridge, 1993,
  pp.~1--295. \MR{95m:20041}

\bibitem[Gui80]{Guichardet(1980)}
A.~Guichardet, \emph{Cohomologie des groupes topologiques et des alg\`ebres de
  {L}ie}, CEDIC, Paris, 1980. \MR{83f:22004}

\bibitem[KR97]{kadison-ringrose(1997a)}
Richard~V. Kadison and John~R. Ringrose, \emph{Fundamentals of the theory of
  operator algebras. {V}ol. {I}}, Graduate Studies in Mathematics, vol.~15,
  American Mathematical Society, Providence, RI, 1997, Elementary theory,
  Reprint of the 1983 original. \MR{98f:46001a}

\bibitem[Lam99]{lam(1999)}
T.~Y. Lam, \emph{Lectures on modules and rings}, Graduate Texts in Mathematics,
  vol. 189, Springer-Verlag, New York, 1999. \MR{99i:16001}

\bibitem[LLS03]{linnell-lueck-schick(2002)}
Peter~A. Linnell, Wolfgang L{\"u}ck, and Thomas Schick, \emph{The {O}re
  condition, affiliated operators, and the lamplighter group}, High-dimensional
  manifold topology, World Sci. Publishing, River Edge, NJ, 2003, pp.~315--321.
  \MR{2 048 726}

\bibitem[LRS99]{Lueck-Reich-Schick(1999)}
Wolfgang L{\"u}ck, Holger Reich, and Thomas Schick, \emph{Novikov-{S}hubin
  invariants for arbitrary group actions and their positivity}, Tel Aviv
  Topology Conference: Rothenberg Festschrift (1998), Amer. Math. Soc.,
  Providence, RI, 1999, pp.~159--176. \MR{2000j:55026}

\bibitem[L{\"u}c97]{Lueck(1997a)}
Wolfgang L{\"u}ck, \emph{Hilbert modules and modules over finite von {N}eumann
  algebras and applications to ${L}\sp 2$-invariants}, Math. Ann. \textbf{309}
  (1997), no.~2, 247--285. \MR{99d:58169}

\bibitem[L{\"u}c98a]{Lueck(1998a)}
\bysame, \emph{Dimension theory of arbitrary modules over finite von {N}eumann
  algebras and ${L}\sp 2$-{B}etti numbers. {I}. {F}oundations}, J. Reine Angew.
  Math. \textbf{495} (1998), 135--162. \MR{99k:58176}

\bibitem[L{\"u}c98b]{Lueck(1998b)}
\bysame, \emph{Dimension theory of arbitrary modules over finite von {N}eumann
  algebras and ${L}\sp 2$-{B}etti numbers. {I}{I}. {A}pplications to
  {G}rothendieck groups, ${L}\sp 2$-{E}uler characteristics and {B}urnside
  groups}, J. Reine Angew. Math. \textbf{496} (1998), 213--236. \MR{99k:58177}

\bibitem[L{\"u}c02]{Lueck(2002)}
\bysame, \emph{{$L\sp 2$}-invariants: theory and applications to geometry and
  {$K$}-theory}, Ergebnisse der Mathematik und ihrer Grenzgebiete. 3. Folge. A
  Series of Modern Surveys in Mathematics [Results in Mathematics and Related
  Areas. 3rd Series. A Series of Modern Surveys in Mathematics], vol.~44,
  Springer-Verlag, Berlin, 2002. \MR{1 926 649}

\bibitem[Mv36]{murray-neumann(1936)}
F.J. Murray and J.~{von Neumann}, \emph{On rings of operators}, Ann. of Math.
  \textbf{37} (1936), 116--229.

\bibitem[Nom54]{nomizu(1954)}
Katsumi Nomizu, \emph{On the cohomology of compact homogeneous spaces of
  nilpotent {L}ie groups}, Ann. of Math. (2) \textbf{59} (1954), 531--538.
  \MR{16,219c}

\bibitem[NS86a]{Novikov-Shubin(1986b)}
S.~P. Novikov and M.~A. Shubin, \emph{Morse inequalities and von {N}eumann
  ${II}\sb 1$-factors}, Dokl. Akad. Nauk SSSR \textbf{289} (1986), no.~2,
  289--292. \MR{88c:58065}

\bibitem[NS86b]{Novikov-Shubin(1986a)}
Sergei~P. Novikov and M.~A. Shubin, \emph{Morse inequalities and von {N}eumann
  invariants of non-simply connected manifolds}, Uspekhi. Matem. Nauk
  \textbf{41} (1986), no.~5, 222--223, in Russian.

\bibitem[Pan89]{pansu(1989)}
Pierre Pansu, \emph{M\'etriques de {C}arnot-{C}arath\'eodory et
  quasiisom\'etries des espaces sym\'etriques de rang un}, Ann. of Math. (2)
  \textbf{129} (1989), no.~1, 1--60. \MR{90e:53058}

\bibitem[Rag72]{raghunathan(1972)}
M.~S. Raghunathan, \emph{Discrete subgroups of {L}ie groups}, Springer-Verlag,
  New York, 1972. \MR{58 \#22394a}

\bibitem[Rei99]{Reich(1999)}
H.~Reich, \emph{Group von neumann algebras and related algebras}, Dissertation,
  Universit\"at G\"ottingen, 1999,
  \url{www.math.uni-muenster.de/u/reichh/publ/diss/diss.dvi}.

\bibitem[Rei01]{Reich(2001)}
Holger Reich, \emph{On the {$K$}- and {$L$}-theory of the algebra of operators
  affiliated to a finite von {N}eumann algebra}, $K$-Theory \textbf{24} (2001),
  no.~4, 303--326. \MR{2003m:46103}

\bibitem[Sha]{shalom(2003)}
Yehuda Shalom, \emph{Harmonic analysis, cohomology and the large scale geometry
  of amenable groups}, preprint, to appear in Acta Math.

\bibitem[Sta70]{stammbach(1970)}
Urs Stammbach, \emph{On the weak homological dimension of the group algebra of
  solvable groups}, J. London Math. Soc. (2) \textbf{2} (1970), 567--570.
  \MR{41 \#8526}

\bibitem[Wei94]{Weibel(1994)}
Charles~A. Weibel, \emph{An introduction to homological algebra}, Cambridge
  University Press, Cambridge, 1994. \MR{95f:18001}

\end{thebibliography}
\end{document}